\documentclass[aos,seceqn,citesort,dvips]{arximspdf}
\usepackage{mathrsfs}
\usepackage{graphicx}

\doi{10.1214/10-AOS820}
\volume{39}
\issue{2}
\pubyear{2011}
\firstpage{956}
\lastpage{989}

\makeatletter

\newtheorem{theorem}{Theorem}[section]
\newtheorem{lem}[theorem]{Lemma}

\newproclaim{rem}{Remark}[section]

\makeatother

\begin{document}
\begin{frontmatter}

\title{A two-stage hybrid procedure for estimating an inverse regression function}
\runtitle{Two-stage estimation of inverse regression function}

\begin{aug}
\author[A]{\fnms{Runlong} \snm{Tang}\ead[label=e1]{rltang@umich.edu}%
\ead[label=u1,url]{http://www.stat.lsa.umich.edu/\textasciitilde rltang}},
\author[A]{\fnms{Moulinath} \snm{Banerjee}\thanksref{t2}\ead[label=e2]{moulib@umich.edu}%
\ead[label=u2,url]{http://www.stat.lsa.umich.edu/\textasciitilde moulib}}
and
\author[A]{\fnms{George} \snm{Michailidis}\corref{}\thanksref{t3}\ead[label=e3]{gmichail@umich.edu}%
\ead[label=u3,url]{http://www.stat.lsa.umich.edu/\textasciitilde gmichail}}
\runauthor{R. Tang, M. Banerjee and G. Michailidis}
\affiliation{University of Michigan}
\address[A]{Department of Statistics\\
University of Michigan\\
1085 South University\\
Ann Arbor, Michigan 48109-1107\\
USA\\
\printead{e1} \\
\phantom{E-mail: }\printead*{e2} \\
\phantom{E-mail: }\printead*{e3} \\
\printead{u1} \\
\phantom{URL: }\printead*{u2} \\
\phantom{URL: }\printead*{u3}} %adresu isvedimo komanda gale!
\end{aug}

\thankstext{t2}{Supported in part by NSF Grant DMS-07-05288.}
\thankstext{t3}{Supported in part by NIH Grant 1RC1CA145444-01.}

% HISTORY:
\received{\smonth{6} \syear{2009}}
\revised{\smonth{2} \syear{2010}}

% ABSTRACT
%
\begin{abstract}
We consider a two-stage procedure (TSP) for estimating an inverse
regression function at a given point, where isotonic regression is used
at stage one to obtain an initial estimate and a local linear
approximation in the vicinity of this estimate is used at stage two. We
establish that the convergence rate of the second-stage estimate can
attain the parametric $n^{1/2}$ rate. Furthermore, a bootstrapped
variant of TSP (BTSP) is introduced and its consistency properties
studied. This variant manages to overcome the slow speed of the
convergence in distribution and the estimation of the derivative of the
regression function at the unknown target quantity. Finally, the finite
sample performance of BTSP is studied through simulations and the
method is illustrated on a data set.
\end{abstract}

% KEYWORDS
%

\end{frontmatter}

%s1 ###
\section{Introduction}\label{sec_intro}

The problem of estimating an \textit{inverse} regression function has a
long history in Statistics,
due to its importance in diverse areas including toxicology,
drug development and engineering. The canonical formulation of the
problem is as follows. Let
\[
Y = f(x) +\varepsilon,
\]
where $f$ is a \textit{monotone} function establishing the relationship
between the design variable $x$ and the response $Y$, and $\varepsilon$
an error term with zero mean and finite variance $\sigma^{2}$.
Further, without loss of generality, it is assumed that $f$ is
isotonic and $x\in[0,1]$. It is of interest to estimate
$d_0=f^{-1}(\theta_0)$ for some $\theta_0$ in the interior of the
range of $f$,
given $f'(d_{0})>0$.

Depending on the nature of the problem, one usually first obtains an
estimate of $f$ and subsequently of $d_{0}$, either from
observational data or from design studies~\cite{Morgan1992}. In the
latter case, one specifies a number of values for the design
variable, and obtains the corresponding responses, which are then
used to get the estimates.

\begin{figure}

\includegraphics{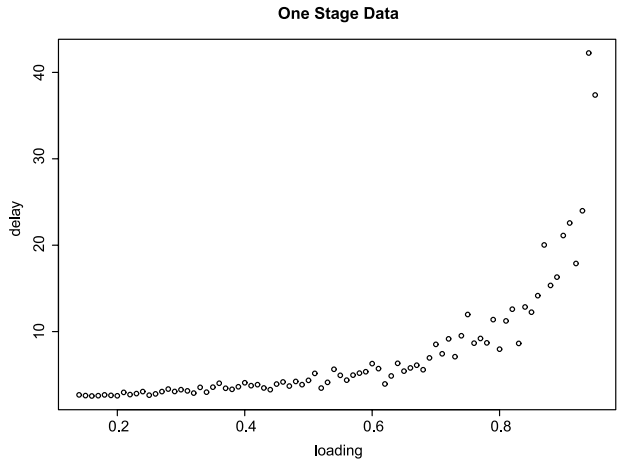}

\caption{The average delay as a function of system's loading.}
\label{Figure Data Analysis Data}
\end{figure}

Motivated by an engineering application, fully described
in Section \ref{data_application}, we introduce a
two-stage design for estimating $d_0$. Specifically, we consider a~complex
queueing system operating in discrete time
under a throughput (average number of customers processed per unit of time)
maximizing resource allocation policy (for details, see Bambos and
Michailidis \cite{Bambos2004}).
Unfortunately the customers' average delay, which is an important
``quality-of-service'' metric of the performance of the system, is not
analytically tractable and can only be obtained via expensive simulations.
The average delay as a function of the
system's loading (number of customers arriving per unit of time) is
depicted in Figure \ref{Figure Data Analysis Data}. The relationship
between system loading and average delay cannot be easily captured
by a simple parametric model; hence, a nonparametric estimator might
be more useful. In addition, given that the responses are obtained
through simulation, only a relatively small number of simulation
runs can be performed. It is of great interest for the system's
operator to obtain accurate estimates of the loading corresponding
to prespecified delay thresholds (e.g., 10 and 15 time units), so as to be
able to decide whether to upgrade the available resources.

%f1 ###
%

The main idea of the proposed two-stage approach is summarized next:
at stage one, an initial set of design
points and their corresponding responses are generated.
Then a first-stage nonparametric estimate of $f$ is obtained and subsequently
a first-stage estimate of $d_{0}$.
Next, a second-stage sampling interval covering $d_{0}$ with high probability
is specified and all new design points are
laid down at the two boundary points of this interval and their
responses obtained.
Finally, a linear regression model
is fitted to the second-stage data by least squares
and a second-stage estimate of $d_0$ computed as the inverse of the
locally approximating line of $f$ at $\theta_{0}$.
As we will see, the employment of a local linear approximation at stage two
allows the second-stage estimate of $d_{0}$ to attain a $\sqrt{n}$
parametric rate of convergence,
despite the nonparametric nature of the problem. To overcome estimation
of several tuning
parameters required by the second-stage estimate, a bootstrapped
variant is introduced and
its consistency properties established. To clinch the asymptotic
results of the proposed
two-stage estimate and its bootstrapped counterpart, a number of subtle
technical issues need
to be addressed and these are resolved in subsequent sections. Before
proceeding further, it
is important to draw attention to the fact that our proposed two-stage
method relies critically on the
\textit{reproducibility} of the experiment: that is, at any stage, it is
possible to sample responses from any
pre-specified covariate value. While reproducibility in this sense is
guaranteed for our motivating application,
the two-stage procedure above is not applicable in the absence of
adequate degree of control on the covariate.
For example, if the covariate is time, the implementation
of a two-stage procedure would require one to go back and sample from
the past, which is impossible.

Isotonic regression is a conceptually natural
and computationally efficient
estimation method for shape-restricted problems
\cite{Barlow1972,Robertson1988}. In the framework of regression,
the asymptotic distribution for the isotonic regression estimator at
a fixed point was first derived in Brunk \cite{Brunk1970},
and then extended in Wright \cite{Wright1981}
and Leurgans \cite{Leurgans1982}. The asymptotic distribution for
the $\mathbb{L}_{1}$-distance between the isotonic estimator and the
regression function was obtained in Durot \cite{Durot2002},
paralleling Groeneboom, Hooghiemstra and Lopuha\"a \cite{Groeneboom1999a} on a~unimodal
density, and then extended in Durot \cite{Durot2007,Durot2008}.
Banerjee and Wellner \cite{Banerjee2005} derive the asymptotic
distribution for the inverse of the distribution function of the
survival time at a given point in the current status model; the
regression version of this result will be used to derive the
asymptotics for the two stage procedures.

The inverse regression problem has been extensively studied in the
context of different applications. For example, in statistical
calibration, the goal is to estimate a scalar quantity $d_0$ from a
model $Z=f(d_0)+\varepsilon$, with $Z$ observed. The information about
the underlying function $f$ comes from experimental data $\{Y_i,
X_i\}$ that follow the same regression model; namely,
$Y_i=f(X_i)+\varepsilon_i$. Osborne \cite{Osborne1991} gives a
comprehensive review of this topic and Gruet \cite{Gruet1996}
provides a kernel based direct nonparametric estimator of $d_{0}$.
It is clear that when $\varepsilon=0$, the calibration problem becomes
the canonical problem described above.

Another active area is provided by the model-based dose-finding
problems in
toxicology and drug development, where $d_0$ corresponds to either
the maximal tolerated dose or the effective dose with respect to a
given maximal toxicity or an efficacy level. Possible dose levels
are often prespecified. The dose-response relationship is usually
assumed to be monotone and described either by parametric models
(e.g., probit, logit \cite{Morgan1992}, multihit \cite{Rai1981},
cubic logistic \cite{Morgan1985}), or by nonparametric models,
for which kernel estimates \cite{Staniswalis1988} and isotonic
regression \cite{Stylianou2002} are employed. On the other hand, due
to ethical and budget considerations, most studies encompass
sequential designs, so that more subjects (e.g., patients) receive
doses close to the target $d_{0}$ (see Rosenberger
\cite{Rosenberger1996} and Rosenberger and Haines
\cite{Rosenberger2002} for comprehensive reviews on the subject).
Stylianou and Flournoy \cite{Stylianou2002} compare parametric
estimators using maximum likelihood and weighted least squares based
on the logit model and nonparametric ones using sample mean and
isotonic regression with a sequential up-and-down biased coin
design, and show that a linearly interpolated isotonic regression
estimator performs best in most simulated scenarios. Further,
Ivanova et al. \cite{Ivanova2003a} claim that the isotonic
regression estimator still performs best for small to moderate
sample sizes with several sequential designs from a family of
up-and-down designs; Gezmu and Flournoy \cite{Gezmu2006} recommend
using smoothed isotonic regression with their group up-and down
designs. All these partially motivate the use of isotonic
regression in our two-stage procedure, though it should be noted
that our approach is markedly different from the ones discussed
above, owing to the different nature of the motivating application;
in particular, ethical constraints that prevent administration of
high dose-levels are absent in our situation.

In a nonparametric setting, one could also employ a fully sequential
Robbins--Monro procedure \cite{Robbins1951} for finding $d_0$. This
corresponds to a stochastic version of Newton's scheme for root
finding problems. Anbar \cite{Anbar1977} considered a modified
Robbins--Monro type procedure approximating the root from one side. A
good review of this area is provided in Lai \cite{Lai2003a}, in
which it is also pointed out that the procedure usually exhibits an
``unsatisfactory finite-sample performance except for linear
problems'' especially in noisy settings, due to the fact that it does
not incorporate modeling for (re)using all the available---up to
that instance---data. Another downside of a sequential design, as
opposed to the \textit{batch} design employed in this study, is the
time and effort required to collect the data points
\cite{muller1990}.

The remainder of the paper is organized as follows: Section \ref{Two
Procedures}
describes the two-stage procedures. The asymptotic properties of the
two-stage estimators are obtained in Section \ref{Sec_Asymptotic_Properties}.
Simulation studies and
data analysis are presented in Sections~\ref{Numerical Work and Applications}
and \ref{data_application}, respectively. We close with
a discussion in Section \ref{Sec_conclusion},
which is followed by an \hyperref[app]{Appendix} containing technical details.

%%%%%%%%%%%%%%%%%%%%%%%%%%%%%%%%%%%%%%%%%%%%%%%%%%%%%%%%%%%%
%%%\input{Adaptive_Procedures}
%s2 ###
\section{Two-stage procedures}\label{Two Procedures}

In this section, we review some necessary background material and
introduce the proposed two-stage estimation procedures.

%s2.1 ###
\subsection{Preliminaries: A single-stage procedure}\label{A Single
Stage Procedure}
We review some material on estimating the parameter of interest
$d_0$ by using isotonic regression combined with a single-stage
design. The procedure is outlined next:
\begin{enumerate}
\item
Choose $n$ increasing design points $\{x_{in}\}_{i=1}^{n}\in[0,1]$ and obtain
the corresponding responses that are generated according to
$Y_{in}=f(x_{in})+\varepsilon_{in}, i=1,2,\ldots,n$, where $f$ is in
$\mathscr{F}_{0}$, a class of increasing real functions on $[0,1]$
with positive and continuous first derivatives in a neighborhood of
$d_{0}$ and $\varepsilon_{in}$ are independently and identically
distributed (i.i.d.) random errors with mean zero and constant variance
$\sigma^2$. Note that the subscript $n$ will be suppressed from now
on for simplicity of notation.
\item
Obtain the isotonic regression estimate $\hat f$ of $f$ from the
data $\{(x_{i},Y_{i})\}_{i=1}^{n}$. (For details, see, Chapter 1 of
Robertson, Wright and Dykstra \cite{Robertson1988}.)
\item
Estimate $d_{0}$ by $\hat d_{n}^{(1)}=\hat f^{-1}(\theta_{0})
=\inf\{x\in[0,1]\dvtx\hat f(x)\geq\theta_{0}\}$, where $\theta_{0}=f(d_{0})$.
\end{enumerate}

In order to study the properties of $\hat f$ and $\hat d_{n}^{(1)}$,
we consider the following further assumption on the design points.
\begin{enumerate}[(A1)]
\item[(A1)] There exists a distribution function $G$,
whose Lebesgue density $g$ is positive at $d_{0}$, such that
$
{\sup_{x\in[0,1]}}|F_{n}(x)-G(x)|=o(n^{-1/3})
$,
where $F_{n}$ is the empirical function of $\{x_{i}\}_{i=1}^{n}$.
\end{enumerate}
For example, the discrete uniform design $x_{i}=i/n$
for $i=1, 2,\ldots, n$ satisfies (A1) with $G$ being
the uniform distribution on $[0,1]$
and $g(d_{0})=1>0$.

The following basic result provides the asymptotic distribution of
$\hat d_{n}^{(1)}$.
\begin{theorem}\label{Thm_one_stage_IR_estimator_d0}
If $f\in\mathscr{F}_{0}$ and \textup{(A1)} holds,
\[
n^{1/3}\bigl(\hat d_{n}^{(1)}-d_{0}\bigr)
\stackrel{d}{\rightarrow}C\mathbb{Z},
\]
where $C=[4\sigma^{2}/(f'(d_{0})^{2}g(d_{0}))]^{1/3}$ and
$\mathbb{Z}$ follows Chernoff's distribution.
\end{theorem}
\begin{rem}
Chernoff's distribution is the distribution of the almost sure unique
maximizer of $B(t)-t^{2}$ on $\mathbb{R}$,
where $B(t)$ denotes a two-sided standard Brownian motion starting at
the origin [$B(0)=0$].
It is symmetric around zero, with tails dwindling faster than those of
the Gaussian and its quantiles
have been tabled in Groeneboom and Wellner \cite{Groeneboom2001}.
\end{rem}

The proof of Theorem \ref{Thm_one_stage_IR_estimator_d0}
follows by adaptations of the arguments from Theorem 1 in Banerjee and Wellner
\cite{Banerjee2005} to the current regression setting. Hence, an approximate
confidence interval for $d_0$ with significance level $1-2\alpha$ can
be constructed
as follows:
%
%e2.1 ###
%
\begin{equation}\label{defn_CI_one_stage_isotonic_regression}
\bigl[\hat d_{n}^{(1)}-n^{-1/3}\hat Cq_{\alpha},
\hat d_{n}^{(1)}+n^{-1/3}\hat Cq_{\alpha}\bigr]
\cap(0,1),
\end{equation}
where $q_\alpha$ denotes the upper $\alpha$ quantile of Chernoff's
distribution and $\hat C$ is a consistent estimate of $C$.

In the presence of relatively small budgets for design points, the
slow convergence rate and the need to estimate $f'(d_0)$ adversely
impact the performance of this procedure. In order to accelerate the
convergence rate, we propose next an alternative that is based on a
two-stage sampling design and uses local linear approximation for
$f$ in stage two.

%s2.2 ###
\subsection{Procedures based on two-stage sampling designs}
\label{Sec_Two_Stage_Plans}
We describe next a~hybrid estimation procedure for estimating $d_0$
based on a two-stage sampling design.
Suppose that the total budget consists of $n$ doses that are going to be
allocated in two stages.
\begin{enumerate}
\item
Allocate\vspace*{1pt} $n_1=np,p\in(0,1)$ design points and obtain the
first-stage data $\{(x_{i},Y_{i})\}_{i=1}^{n_{1}}$, the isotonic
regression estimate of $f$ and the estimate $\hat d_{n_{1}}^{(1)}$
of $d_0$ as outlined in Section \ref{A Single Stage Procedure}. Note
that by $np$, we denote by
$\lfloor np \rfloor$ or $\lfloor np \rfloor+ 1$, depending on whether
$n-\lfloor np \rfloor$ is
even or not. Also, recall that the additional subscript $n$ is suppressed.
\item
Determine two second-stage
sampling points $L$ and $U$ symmetrically around $\hat d_{n_1}^{(1)}$,
where $L=\hat d_{n_{1}}^{(1)}- Kn_{1}^{-\gamma}$ and
$U=\hat d_{n_{1}}^{(1)}+ Kn_{1}^{-\gamma}$, for some constants
$\gamma>0$ and $K>0$.
\item
Allocate the remaining $n-n_1$ design points
\textit{equally} to $L$ and $U$ and
generate the responses as $Y'_{i}=f(L)+\varepsilon'_{i}$ and
$Y''_{i}=f(U)+\varepsilon''_{i}$ for $i=1,2,\ldots,n_{2}$, with
$\{\varepsilon'_i\}$ and $\{\varepsilon''_i\}$ being i.i.d. random errors
with mean zero and constant variance~$\sigma^2$, mutually independent
and also independent of~$\{\varepsilon_i\}$.
\item
Fit the second-stage data $\{(L,Y'_{i}),(U,Y''_{i})\}$
with the linear model $y=\beta_{0}+\beta_{1}x$
using least\vspace*{1pt} squares.
Denote the resulting intercept and slope estimates by
$(\hat\beta_{0}, \hat\beta_{1})$, respectively. Then,
the second-stage (or two-stage) estimator of $d_{0}$ is given
by $\tilde d_{n}^{(2)}=(\theta_{0}-\hat\beta_{0})/\hat
\beta_{1}$.
\end{enumerate}

Asymptotic properties of $\tilde d_{n}^{(2)}$ will be established in
Section \ref{Two-Stage Estimator}. For example,
when $f$ is in a subset of $\mathscr{F}_{0}$, denoted as
$\mathscr{F}$, the third derivatives of whose elements are uniformly
bounded around $d_{0}$,
and $\gamma\in(1/4,1/3)$, we have
%
%e2.2 ###
%
\begin{equation}\label{Result tilde dn limit distribution}
n^{1/2}\bigl(\tilde
d_{n}^{(2)}-d_{0}\bigr)\stackrel{d}{\rightarrow}
\frac{\sigma}{f'(d_{0})(1-p)^{1/2}}N(0,1),
\end{equation}
where $\stackrel{d}{\rightarrow}$ denotes convergence in distribution.
Thus, the convergence rate of the two-stage estimator of $d_{0}$
becomes $n^{1/2}$,
the standard parametric convergence rate,
which is faster than the $n^{1/3}$ convergence rate of the one-stage
isotonic regression estimator.

However, when constructing confidence intervals from asymptotic results
like~(\ref{Result tilde dn limit distribution}),
we face two difficulties. One is that
the limiting distributions of interest still depend on $f'(d_{0})$,
accurate estimation of which is difficult for small to moderate sample sizes.
The other one, which is less obvious but perhaps with more serious
practical implications, is that
the asymptotic results of interest suffer slow speed of convergence in
distribution.
Therefore, a bootstrap variant of the two-stage procedure that avoids
direct estimation of $f'(d_0)$ is introduced
and is seen to relieve the slow convergence problem.

%s2.3 ###
\subsection{Bootstrapping the two-stage estimator}
\label{Bootstrapping The Two-Stage Estimator}
The steps of the bootstrapped two-stage procedure are outlined next.
\begin{enumerate}
\item
Follow steps 1--4 to obtain the second stage design points $L$ and $U$,
responses
$\{Y'_i\}$ and $\{Y^{\prime\prime}_i\}$ and $\tilde d_{n}^{(2)}$.
\item
Sample with replacement, responses $\{Y_{i}^{\prime\star}\}_{i=1}^{n_{2}}$
and $\{Y_{i}^{\prime\prime\star}\}_{i=1}^{n_{2}}$,
from $\{Y'_{i}\}_{i=1}^{n_{2}}$ and $\{Y''_{i}\}_{i=1}^{n_{2}}$, respectively.
Construct the corresponding bootstrapped second-stage (or two-stage)
estimator $\tilde d_{n}^{(2)\star}$,
and calculate the corresponding root $R_{n}^{\star}=n^{1/2}(\tilde
d_{n}^{(2)\star}-\tilde d_{n}^{(2)})$.
\item
Repeat the previous step $B$ times to obtain $\{R_{n}^{\star b}\}
_{b=1}^{B}$. Subsequently, calculate
the lower and upper $\alpha$ quantiles, $q_{l}^{\star}$ and
$q_{u}^{\star}$,
of $\{R_{n}^{\star b}\}_{b=1}^{B}$. Finally, construct a $1-2\alpha$
bootstrapped Wald-type
confidence interval for $d_{0}$ as
%
%e2.3 ###
%
\begin{equation}\label{BLTSP_2nd_CI}
\bigl[\tilde d_{n}^{(2)}-n^{-1/2}q_{u}^{\star}, \tilde
d_{n}^{(2)}-n^{-1/2}q_{l}^{\star}\bigr].
\end{equation}
\end{enumerate}
Note that the procedure does not require estimation of $f'(d_0)$.

The asymptotic properties of the bootstrapped two-stage estimator are
established in Section
\ref{Bootstrapped Two-Stage Estimator}.
For example, when $f\in\mathscr{F}$, $\gamma\in(0,1/3)$ and
all the absolute moments of the random error are finite,
we have
%
%e2.4 ###
%
\begin{equation}\label{Result_tilde_dn_star_limit_distribution}
n^{1/2}\bigl(\tilde
d_{n}^{(2)\star}-\tilde d_{n}^{(2)}\bigr)\stackrel{d^{\star}}{\rightarrow}
\frac{\sigma}{f'(d_{0})(1-p)^{1/2}}N(0,1),\qquad
(P\mbox{-a.s.}),
\end{equation}
where $\stackrel{d^{\star}}{\rightarrow}$ implies convergence in distribution
conditional on the data obtained from the employed two-stage design.

From (\ref{Result tilde dn limit distribution})
and (\ref{Result_tilde_dn_star_limit_distribution}),
the strong consistency of the bootstrapped estimator $\tilde
d_{n}^{(2)\star}$ is ensured
for $f\in\mathscr{F}$ and $\gamma\in(1/4,1/3)$.
In fact,
the strong assumption on the random error
can be replaced by a mild one that the sixth moment of
the random error is finite,
at the price of replacing strong consistency with weak consistency.
Therefore, the bootstrapped procedure is theoretically validated under
certain conditions.
\begin{rem}
Both the two-stage estimator and its bootstrapped variant
rely on the choice of a number of tuning parameters: $p$, $\gamma$ and $K$.
Practical procedures for their selection will be discussed in Section
\ref{Numerical Work and Applications}.
\end{rem}

%%%\input{Asymptotic_Results}
%s3 ###
\section{Asymptotic properties of two-stage estimators}
\label{Sec_Asymptotic_Properties}

In this section, we establish the asymptotic properties of both the
two-stage estimator and its bootstrapped variant for $d_0$. We start
by discussing the two-stage estimator~$\tilde{d}^{(2)}_n$.

%s3.1 ###
\subsection{Two-stage estimator}\label{Two-Stage Estimator}
All results in this subsection are derived under the assumption
(A1). According to the two-stage procedure,
\[
(\hat\beta_{0}, \hat\beta_{1})=\mathop{\arg\min}_{\beta
_{0},\beta_{1}\in\mathbb{R}}
\sum_{i=1}^{n_{2}}[(Y'_{i}-\beta_{0}-\beta_{1}L)^{2}+(Y''_{i}-\beta
_{0}-\beta_{1}U)^{2}].
\]
Denote $Y_{i}^{+}=Y''_{i}+Y'_{i}$ and $Y_{i}^{-}=Y''_{i}-Y'_{i}$.
Then,
%
%e3.1 ###
%
\begin{equation}\label{TED G beta0 beta1} %Two End Design General beta0
\hat\beta_{0}
= (2n_{2})^{-1}\sum_{i=1}^{n_{2}}Y_{i}^{+}-\hat d_{n_{1}}^{(1)} \hat
\beta_{1},\qquad
\hat\beta_{1}=(2Kn_{1}^{-\gamma}n_{2})^{-1}\sum_{i=1}^{n_{2}}Y_{i}^{-}.
\end{equation}
Setting $\theta_{0}=\hat\beta_{0}+ \hat\beta_{1}\tilde
d_{n}^{(2)}$ gives
%
%e3.2 ###
%
\begin{equation}\label{TED_G_tilded}
\tilde d_{n}^{(2)} = (1/\hat\beta_{1})(\theta_{0}-\hat\beta_{0})
= (1/\hat\beta_{1})\Biggl[\theta_{0}-(2n_{2})^{-1}\sum_{i=1}^{n_{2}}Y_{i}^{+}\Biggr]
+ \hat d_{n_{1}}^{(1)}.
\end{equation}

In order to analyze $\tilde d_{n}^{(2)}$, additional assumptions
about the smoothness of the underlying function $f$ around $d_0$ are
required. We consider the following three classes of underlying
functions:
\begin{eqnarray*}
\mathscr{F}
&=& \{f\in\mathscr{F}_{0}\dvtx f'''(x) \mbox{ is } \mathrm{UBN}(d_{0})\}, \\
\mathscr{F}_{1}
&=& \{f\in\mathscr{F}_{0}\dvtx f''(d_{0})\not= 0, f'''(x) \mbox{ is }
\mathrm{UBN}(d_{0})\}, \\
\mathscr{F}_{2}
&=& \bigl\{f\in\mathscr{F}_{0}\dvtx f''(d_{0}) = 0, f'''(d_{0})\not= 0,
f^{(4)}(x)\mbox{ is } \mathrm{UBN}(d_{0})\bigr\},
\end{eqnarray*}
where $\operatorname{UBN}(d_{0})$ means ``uniformly bounded in a neighborhood of
$d_{0}$.''
Then, the mutually exclusive $\mathscr{F}_{1}$ and $\mathscr{F}_{2}$ are
subsets of $\mathscr{F}$.
\begin{rem}
A function in $\mathscr{F}_{2}$ is exactly locally linear at $d_{0}$
while that in $\mathscr{F}_{1}$ is not. Notice that both
$\mathscr{F}_{2}$ and $\mathscr{F}_{1}$ depend on $d_{0}$. For
example, consider the sigmoid function
$f(x)=\exp\{a(x-b)\}/(1+\exp\{a(x-b)\})$ for some constants $a>0$ and
$b\in(0,1)$. It belongs to $\mathscr{F}_{2}$ if $d_{0}=b$ and to
$\mathscr{F}_{1}$ otherwise. Obviously, the size of
$\mathscr{F}_{2}$ is much smaller than that of $\mathscr{F}_{1}$.
However, the asymptotic results for $f\in\mathscr{F}_{2}$ should also
provide good approximations for functions that are approximately
linear in the vicinity of $d_{0}$. Hence, the class
$\mathscr{F}_{2}$ is also of interest.
\end{rem}

We consider next the asymptotic properties of $\tilde d_{n}^{(2)}$,
starting with the consistency of the two-stage estimator.
\begin{lem}\label{Two_Ends_General_Case_One_Lemma_Consistency}
For $f \in\mathscr{F}$ and $\gamma\in(0,1/2)$, we have
\[
\hat\beta_{0}\stackrel{P}{\rightarrow} f(d_{0})-f'(d_{0})d_{0},\qquad
\hat\beta_{1}\stackrel{P}{\rightarrow} f'(d_{0})\quad \mbox{and}\quad
\tilde d_{n}^{(2)} \stackrel{P}{\rightarrow} d_{0}.
\]
\end{lem}

Based on Lemma \ref{Two_Ends_General_Case_One_Lemma_Consistency}, we obtain
the asymptotic distribution of $\tilde{d}^{(2)}_n$ in the next theorem.
It turns out that the asymptotic results with $f \in\mathscr{F}_{1}$ and
$\mathscr{F}_{2}$ are the same for $\gamma>1/4$.
This implies that
the nonlinearity of $f$ at $d_{0}$ becomes asymptotically ignorable
as the length of the neighborhood of $d_{0}$
shrinks fast enough.
\begin{theorem}\label{thm two stage estimator case one LD}
For $f\in\mathscr{F}$ and $\gamma\in(1/4, 1/2)$,
\begin{eqnarray*}
n^{1/2}\bigl(\tilde d_{n}^{(2)}-d_{0}\bigr) &\stackrel{d}{\rightarrow}&
C_{2}Z_{1}\qquad\mbox{for } \gamma\in(1/4, 1/3), \\
n^{1/2}\bigl(\tilde d_{n}^{(2)}-d_{0}\bigr) &\stackrel{d}{\rightarrow}&
C_{2}Z_{1}+C_{3}\mathbb{Z}Z_{2}\qquad \mbox{for } \gamma=1/3, \\
n^{(5/6 -\gamma)}\bigl(\tilde d_{n}^{(2)}-d_{0}\bigr) &\stackrel{d}{\rightarrow}&
C_{3}\mathbb{Z}Z_{2}\qquad \mbox{for } \gamma\in(1/3, 1/2);
\end{eqnarray*}
for $f\in\mathscr{F}_{1}$ and $\gamma\in(0, 1/4]$,
\begin{eqnarray*}
n^{2\gamma}\bigl(\tilde d_{n}^{(2)}-d_{0}\bigr) &\stackrel{d}{\rightarrow}&
C_{1}\qquad \mbox{for } \gamma\in(0, 1/4), \\
n^{1/2}\bigl(\tilde d_{n}^{(2)}-d_{0}\bigr) &\stackrel{d}{\rightarrow}&
C_{1}+C_{2}Z_{1}\qquad \mbox{for } \gamma= 1/4;
\end{eqnarray*}
for $f\in\mathscr{F}_{2}$ and $\gamma\in(1/8, 1/4]$,
\[
n^{1/2}\bigl(\tilde d_{n}^{(2)}-d_{0}\bigr) \stackrel{d}{\rightarrow}
C_{2}Z_{1}\qquad \mbox{for } \gamma\in(1/8, 1/4];
\]
where $C_{1}=-K^{2}p^{-2\gamma}f''(d_{0})/[2f'(d_{0})]$,
$C_{2}=\sigma/[f'(d_{0})(1-p)^{1/2}]$, $C_{3}=CC_{2}/\break K$,
$C$ is as given in Theorem \ref{Thm_one_stage_IR_estimator_d0},
$Z_1$ and $Z_2$ are standard normal, $\mathbb{Z}$ follows Chernoff's
distribution
and $\mathbb{Z}, Z_1, Z_2$ are mutually independent.
\end{theorem}
\begin{rem}
Theorem \ref{thm two stage estimator case one LD} characterizes the
convergence rate of the estimator
in terms of the size of the shrinking neighborhood. It shows that
for $\gamma\in[1/4,1/3]$ the parametric rate of $n^{1/2}$ is
achieved given $f\in\mathscr{F}$. On the other hand, for the boundary
values of $\gamma=1/4$
and $1/3$, there exists asymptotic bias in the former case
(for $f\in\mathscr{F}_{1}$), while in
the latter case the asymptotic variance increases. For $\gamma>1/3$, the
rate of convergence falls below $\sqrt{n}$, while for $\gamma<1/4$
and $f\in\mathscr{F}_1$ the limit distribution of the two-stage
estimate is
degenerate and thus not conducive to inference.
Hence, these results suggest selecting $\gamma$ in the $(1/4,1/3)$ range.
Note that,
the function class $\mathscr{F}_{2}$ achieves the $n^{1/2}$ rate of convergence
for a slightly larger range of values for $\gamma$ than $\mathscr{F}_{1}$.
This is a consequence
of the near linearity of $f$ in the vicinity of $d_0$, which allows a
good linear
approximation of $f$ with a relatively long interval $[L,U]$.
\end{rem}
\begin{rem}
The case of $\gamma<1/8$ has been omitted for $f\in\mathscr{F}_{2}$,
since it involves a
Taylor expansion of $f$ up to its fifth derivative. Nevertheless, in
principle no other technical challenges are in play.
\end{rem}

%s3.2 ###
\subsection{Bootstrapped two-stage estimator}\label{Bootstrapped
Two-Stage Estimator}

We consider next the asymptotic properties of the bootstrapped
two-stage estimator,
which is
%
%e3.3 ###
%
\begin{equation}\label{Bootstrap_2nd_stage_estimator_d0}\quad
\tilde d_{n}^{(2)\star}
= (1/\hat\beta_{1}^{\star})(\theta_{0}-\hat\beta_{0}^{\star})
=(1/\hat\beta_{1}^{\star})
\Biggl[f(d_{0})-(2n_{2})^{-1}\sum_{i=1}^{n_{2}}
Y_{i}^{\star+}\Biggr] + \hat d_{n_{1}}^{(1)},
\end{equation}
where
$Y_{i}^{\star+}=Y''^{\star}_{i}+Y'^{\star}_{i}$,
$Y_{i}^{\star-}=Y''^{\star}_{i}-Y'^{\star}_{i}$ and
%
%e3.4 ###
%
\begin{equation}\label{Bootstrap 2nd stage estimator beta 0 beta1}
\hat\beta_{0}^{\star} = (2n_{2})^{-1}\sum_{i=1}^{n_{2}}
Y_{i}^{\star+} - \hat d_{n_{1}}^{(1)}\hat\beta_{1}^{\star},\qquad
\hat\beta_{1}^{\star} = (2Kn_{1}^{-\gamma}n_{2})^{-1}\sum_{i=1}^{n_{2}}
Y_{i}^{\star-}.
\end{equation}

We now present a probabilistic framework needed to clearly
establish the asymptotic properties of the bootstrapped estimator
rigorously. The point is that the design points and random errors
involved in the sampling mechanism are assumed to come from
triangular arrays but not necessarily from sequences.

Let $\{\{ x_{in}\}_{i=1}^{n} \}_{n=1}^\infty$ be a triangluar array of
distinct design
points in $[0,1]$ and $\varepsilon$ a continuous random variable in
$\mathbb{R}$ with mean $0$ and finite variance $\sigma^{2}>0$. Now,
there exists, on some probability space $(\Omega, \mathscr{A}, P)$,
a~set of random errors
$\{ \{\varepsilon_{in}\}_{i=1}^{n},
\{ \varepsilon_{in}^{\prime}\}_{i=1}^{n},
\{\varepsilon_{in}^{\prime\prime}\}_{i=1}^{n}\}_{n=1}^{\infty}$ which are i.i.d.
copies of $\varepsilon$. Then, suppressing the subscript $n$,
$\{\{x_{i}\}_{i=1}^{n_{1}},
\{\varepsilon_{i}(\omega)\}_{i=1}^{n_{1}},
\{\varepsilon_{i}^{\prime}(\omega)\}_{i=1}^{n_{2}},\break
\{\varepsilon_{i}^{\prime\prime}(\omega)\}_{i=1}^{n_{2}}\}_{n=1}^{\infty}$
represents a~fixed triangular array of real numbers for a~fixed
$\omega\in\Omega$, where $n=n_1+2n_2$ with $n_1$ and $2n_2$ denoting
the first and second stage sample sizes.

Given $\omega\in\Omega$, according to the sampling mechanism used
in the bootstrapped procedure, the data obtained from the first
stage are given by $\{(x_{i},Y_{i}(\omega))\}_{i=1}^{n_{1}}$, which
are subsequently used to obtain $\hat d_{n_{1}}^{(1)}(\omega)$ and
the lower and upper boundary points $L(\omega)$ and $U(\omega)$ to
be used in the second stage. Hence, the second-stage data are given
by $\{L(\omega), Y_{i}^{\prime}(\omega)\}$ and
$\{U(\omega),Y_{i}^{\prime\prime}(\omega)\}$ and the resulting estimate by
$\tilde d_{n}^{(2)}(\omega)$. The procedure then requires
bootstrapping $\{Y_{i}^{\prime}(\omega)\}_{i=1}^{n_{2}}$ and
$\{Y_{i}^{\prime\prime}(\omega)\}_{i=1}^{n_{2}}$, which is conceptually
equivalent to bootstrapping
$\{\varepsilon_{i}^{\prime}(\omega)\}_{i=1}^{n_{2}}$ and
$\{\varepsilon_{i}^{\prime\prime}(\omega)\}_{i=1}^{n_{2}}$ to get
$\{\varepsilon_{i}^{\prime\star}\}_{i=1}^{n_{2}}$ and
$\{\varepsilon_{i}^{\prime\prime\star}\}_{i=1}^{n_{2}}$, so that
$Y_{i}^{\prime\star}=f(L(\omega))+\varepsilon_{i}^{\prime\star}$ and
$Y_{i}^{\prime\prime\star}=f(U(\omega))+\varepsilon_{i}^{\prime\prime\star}$ for
$i=1,2,\ldots,n_{2}$. Note that given $\omega$ and $n$, the
bootstrapped second-stage random errors
$\{\varepsilon_{i}^{\prime\star}\}_{i=1}^{n_{2}}$ and
$\{\varepsilon_{i}^{\prime\prime\star}\}_{i=1}^{n_{2}}$ are i.i.d. uniform random
variables supported on $\{\varepsilon_{i}^{\prime}(\omega)\}_{i=1}^{n_{2}}$
and $\{\varepsilon_{i}^{\prime\prime}(\omega)\}_{i=1}^{n_{2}}$, respectively.
Finally, the bootstrapped estimate $\tilde d_{n}^{(2)\star}$ is
calculated from $\{(L(\omega), Y_{i}^{\prime\star}),
(U(\omega),Y_{i}^{\prime\prime\star})\}_{i=1}^{n_{2}}$.

Thus, given $\omega$ and with $n$ increasing, the design points and
random errors are sampled as rows from the fixed triangular array.
Then the bootstrapped random errors
$\{\varepsilon_{i}^{\prime\star}\}_{i=1}^{n_{2}}$ and
$\{\varepsilon_{i}^{\prime\prime\star}\}_{i=1}^{n_{2}}$ also form
triangular arrays as $n$ varies.
Given $\omega$ and $n$, the randomness of $\tilde
d_{n}^{(2)\star}$ comes from the bootstrapping step.

Under the above theoretical setting, in order to obtain the strong
consistency of the bootstrapped estimator, we consider the following
strong assumptions on the design points, the regression function
and the random errors.
\begin{enumerate}[(A1)]
\item[(A2)] There exists a distribution function $G$,
whose Lebesgue density $g$ is positive
and has a bounded first derivative on $[0,1]$, such that
$\sup_{x\in[0,1]}|F_{n}(x)-G(x)|\lesssim n^{-1/2}$,
where $F_{n}$ is the empirical function of $\{x_{i}\}_{i=1}^{n}$
and ``$\lesssim$'' denotes that the left-hand side is less than
a constant times the right-hand side.
\item[(A3)] The regression function $f\in\mathscr{F}_{0}$ is
differentiable on $[0,1]$
with\break $\inf_{x\in[0,1]}f'(x)$ and $\sup_{x\in[0,1]}f'(x)$ both
positive and
finite.
\item[(A4)] All the absolute moments of $\varepsilon$ are finite,
that is, $\mathbb{E}|\varepsilon|^{q}<\infty$ for all $q\in\mathbb{N}$.
\end{enumerate}
\begin{rem}
There exist triangular arrays of design points satisfying~(A2).
For example, with $x_{i}=i/n$
for $i=1,2,\ldots,n$ and all $n$,
we have an array of discrete uniform designs on $[0,1]$.
Let $G$ be the uniform distribution function on $[0,1]$.
Then, for this special array
$\sup_{x\in[0,1]}|F_{n}(x)-G(x)|\leq1/n$.
Note that (A2) is stronger than (A1).
A random variable with finite moment generating function
in a small neighborhood of 0 satisfies (A4), such as a normal random
variable. The assumptions (A2) to (A4) are essentially the fixed
design versions of the assumptions for Lemma 1 of Durot
\cite{Durot2008}, a modification of which enables us to identify a
crucial boundary rate for the almost sure convergence of the
isotonic regression estimator of $d_{0}$. This boundary rate plays a
central role in the strong consistency of the bootstrapped
estimator.
\end{rem}

Next, we state results on the strong consistency of $\hat\beta_{1}$
and the conditional weak consistency of $\hat\beta_{1}^{\star}$ and
then on strong consistency of the bootstrapped estimator. Note that
$P^\star$ denotes the probability of the bootstrapped data
conditional on the original data.
\begin{lem}\label{Lemma_Bootstrapped_Root_Consistency}
If $f\in\mathscr{F}$, $\gamma\in(0,1/2)$ and \textup{(A2)} to \textup{(A4)} hold,
\[
\hat\beta_{1} \rightarrow f'(d_{0}),\qquad  (P\mbox{-a.s.}),\qquad
\hat\beta_{1}^{\star} \stackrel{P^{\star}}{\rightarrow}
f'(d_{0}),\qquad  (P\mbox{-a.s.}),
\]
where $\stackrel{P^{\star}}{\rightarrow}$ denotes convergence in probability
conditional on a given $\omega$.
\end{lem}
\begin{theorem}\label{Thm Bootstrapped two stage estimator LD}
If $f\in\mathscr{F}$, $\gamma\in(0,1/3)$ and \textup{(A2)} to \textup{(A4)} hold,
\[
n^{1/2}\bigl(\tilde
d_{n}^{(2)\star}-\tilde d_{n}^{(2)}\bigr)\stackrel{d^{\star}}{\rightarrow}
C_{2}Z_{1},\qquad
(P\mbox{-a.s.}),
\]
where $C_{2}$ and $Z_{1}$
are as in Theorem \ref{thm two stage estimator case one LD}. That is,
\[
\sup_{t\in\mathbb{R} } \bigl| P^{\star}\bigl( n^{1/2}\bigl(\tilde
d_{n}^{(2)\star}-\tilde d_{n}^{(2)}\bigr)\leq t\bigr)-P( C_{2}Z_{1}\leq t
)\bigr| \stackrel{\mathit{a.s.}}{\rightarrow} 0.
\]
\end{theorem}

From the above strong consistency, the corresponding weak
consistency follows under the same set of conditions. However, weak
consistency can be obtained with the following weaker requirement
on the random error:
\begin{enumerate}[(A1)]
\item[(A5)] The sixth moment of $\varepsilon$ is finite,
that is, $\mathbb{E}\varepsilon^{6} <\infty$.
\end{enumerate}
\begin{theorem}\label{Thm Bootstrapped two stage estimator LD Weak}
If $f\in\mathscr{F}$, $\gamma\in(0,1/3)$ and \textup{(A1)} and \textup{(A5)} hold,
for $t\in\mathbb{R}$,
\[
\sup_{t\in\mathbb{R} } \bigl| P^{\star}\bigl( n^{1/2}\bigl(\tilde
d_{n}^{(2)\star}-\tilde d_{n}^{(2)}\bigr)\leq t\bigr)-P( C_{2}Z_{1}\leq t
)\bigr| \stackrel{P}{\rightarrow} 0,
\]
where $C_{2}$ and $Z_{1}$
are as in Theorem \ref{thm two stage estimator case one LD}.
\end{theorem}
\begin{rem}\label{remark for thm bootstrpped two stage estimaotor LD}
Comparing Theorem \ref{Thm Bootstrapped two stage estimator LD}
with Theorem \ref{thm two stage estimator case one LD},
we see that, under the strong assumption (A5) on the random errors,
the bootstrapped estimator
is strongly consistent for $f\in\mathscr{F}$ and $\gamma\in(1/4,1/3)$,
which is exactly the $\gamma$-range of most interest.
Further, if $f$ is locally linear at $d_{0}$, that is,
$f\in\mathscr{F}_{2}$, the strong consistency continues to hold for
$\gamma\in(1/8, 1/4]$. Similar conclusions on weak consistency hold
by comparing Theorem \ref{Thm Bootstrapped two stage estimator LD
Weak}
with Theorem \ref{thm two stage estimator case one LD},
but under the milder assumption (A5) on the random errors.
\end{rem}

%%%\input{Numerical_Application_GM}
%s4 ###
\section{Performance evaluation}\label{Numerical Work and Applications}

In this section, through an extensive simulation study we investigate
the finite sample performance of the one-stage procedure
(henceforth, OSP), the proposed two-stage procedure (TSP)
and its bootstrapped variant (BTSP).

Notice that for practically implementing the OSP, as well as the
two-stage procedures, estimates\vadjust{\goodbreak} of $f'(d_0)$ and $\sigma^2$ need to
be obtained; the resulting procedures are called POSP, PTSP and
PBTSP, respectively (Practical OSP, TSP and BTSP). For $\sigma^2$,
we employ the nonparametric estimator proposed by Gasser, Sroka and Jennen-Steinmetz
\cite{Gasser1986}, which is based on local linear fitting. Suppose
the data $\{(x_{i}, Y_{i})\}_{i=1}^{n}$ are already sorted in
ascending order of $x_{i}$'s. Then, we calculate
\[
S^{2}=(n_{1}-2)^{-1}\sum_{i=2}^{n-1}c_{i}^{2}\tilde\varepsilon_{i}^{2},
\]
where $\tilde\varepsilon_{i}=a_{i}Y_{i-1}+b_{i}Y_{i+1}-Y_{i}$,
$c_{i}^{2}=(a_{i}^{2}+b_{i}^{2}+1)^{-1}$,
$a_{i}=(x_{i+1}-x_{i})/(x_{i+1}-x_{i-1})$ and
$b_{i}=(x_{i}-x_{i-1})/(x_{i+1}-x_{i-1})$, for $i = 2, 3,\ldots,
n-1$. An estimate of $f'(d_0)$ is obtained through the local
quadratic regression estimator proposed by Fan and Gijbels
\cite{Fan1996}, at the estimate $\hat{d}_{n}^{(1)}$. Specifically,
let $K(\cdot)$ denote the Epanechnikov kernel function and $h>0$ the
bandwidth, so that $K_{h}(\cdot)= (1/h)K(\cdot/h)$. Further, let
$\hat\eta= (\hat\eta_{0}, \hat\eta_{1},\hat\eta_{2})$ be given
by
\[
\hat\eta= \mathop{\arg\min}_{\eta\in\mathbb{R}^{3}}
\sum_{i=1}^{n}\Biggl[ Y_{i} - \sum_{j=0}^{2}\eta_{j}\bigl(x_{i}-\hat
d_{n}^{(1)}\bigr)^{j}\Biggr]^{2}
K_{h}\bigl(x_{i}-\hat d_{n}^{(1)}\bigr).
\]
Then, the local quadratic regression estimator of $f'(\hat
d_{n}^{(1)})$ is given by $\hat\eta_{1}$. The bandwidth $h$ is
chosen by\vspace*{2pt} first fitting a fifth order polynomial function to the
data to obtain $\hat{f}_{\mathrm{pol}}(x)=\sum_{j=0}^{5}\hat\alpha_{j}x^{j}$
. Next, the estimate of the third order derivative of $f$ at
$\hat{d}_{n}^{(1)}$ is obtained by $\hat f^{(3)}_{\mathrm{pol}}(\hat
d_{n}^{(1)})= 6\hat\alpha_{3}+24\hat\alpha_{4}\hat d_{n}^{(1)}
+60\hat\alpha_{5} (\hat d_{n}^{(1)})^{2}$. Finally, the bandwidth
$h$ is calculated as
\[
\hat h_{\mathrm{opt}}= C_{1,2}(K)
\bigl[ S^{2}/\bigl(\hat f^{(3)}_{\mathrm{pol}}\bigl(\hat d_{n}^{(1)}\bigr)\bigr)^{2}
\bigr]^{1/7}n^{-1/7},
\]
where $C_{1,2}(K)=2.275$.

For the two-stage procedures, the tuning parameters
$\gamma$ and $K$ need to be specified for obtaining the
second-stage sampling points $L$ and $U$. We select them as the end
points of a \textit{high level Wald-type confidence interval} calculated
from the
first-stage data; that is, $\gamma$ and $K$ satisfy
%
%e4.1 ###
%
\begin{equation}\label{K gamma CI}
Kn_{1}^{-\gamma} = Cq_{\beta}n_{1}^{-1/3},
\end{equation}
where $q_{\beta}$ is the upper $\beta$ quantile of $\mathbb{Z}$. On
the other hand, a good quantitative rule for selecting the
first-stage sample proportion $p$ is not available; nevertheless, a
practical qualitative rule of thumb dictates that $p$ should
decrease, while $np$ should increase as the sample size increases.
In our simulation study, a number of different values for $p$ are
considered.

Finally, due to presence of small sample sizes the following modification
of the second-stage estimator is adopted:
\[
\tilde d_{n}^{(2)} = \cases{
\min\bigl(\max\bigl((\theta_{0}-\hat\beta_{0})/\hat\beta_{1}, 0\bigr),1\bigr),
&\quad
if $\hat\beta_{1}>0$,\cr
\hat d_{n_{1}}^{(1)}, &\quad otherwise.}
\]
The same modification applies to the bootstrapped second-stage estimator
in BTSP.
\begin{rem}
Note that our method for choosing the tuning parameters $\gamma, K$
brings in
another subjective parameter $\beta$. However, the choice of $\beta$
is guided by a rational
principle, namely the requirement that the chosen interval contain the
truth with high probability.
The magnitude of $\beta$ is related to how conservative the
experimenter wants to be in the construction
of the second stage interval which is fundamentally subjective. Also,
our rule of thumb regarding the choice
of $p$ is based on the idea that with larger budgets smaller $p$'s at
stage one will still lead to reasonably precise sampling
intervals at stage two, leaving a larger proportion of points for stage
two and the possibility of more accurate conclusions.
\end{rem}

The basic settings of the simulation study are as follows: two
regression functions are considered, $f_{1}(x)=x^{2}+x/5$ and
$f_{2}(x)=e^{4(x-0.5)}/(1+e^{4(x-0.5)})$ for $x\in[0,1]$. The
first-stage design points are drawn from a discrete uniform
distribution on $[0,1]$, that is, $x_{i}=i/(n_{1}+1)$. Further, the
target is set to $d_{0}=0.5$, the standard deviation of the random
error $\sigma$ to $0.1$, $0.3$ and $0.5$, the sample size $n$ ranges
from 50 to 500 in increments of 50, while the first-stage sample
proportion $p$ ranges from 0.2 to 0.8 in increments of 0.1. Finally,
the levels of significance $\alpha$ and $\beta$ are set to $0.025$.
Note that $\beta$ is only required to be small and the specific
choice of $0.025$ is somewhat arbitrary. The following quantities are
computed: coverage rates and average lengths of confidence
intervals, and mean squared errors of estimators. The
simulation programs and more results can be found on the first
author's webpage:
\href{http://www.stat.lsa.umich.edu/\textasciitilde rltang}{www.stat.lsa.umich.edu/\textasciitilde rltang}.
In this paper, we show part of the results for saving space.
\begin{rem}
Choosing $\gamma$ and $K$ via equation (\ref{K gamma CI}) is
theoretically equivalent to having $\gamma=1/3$ and $K=Cq_{\beta}$.
Notice that strictly speaking,
neither strong nor weak consistency for $\gamma=1/3$ is
expected to hold for the bootstrapped estimator.
However, it is reasonable to expect that for realistic sample sizes,
the performance of the bootstrap would be satisfactory,
since $\gamma=1/3$ is at the boundary of consistency. The
obtained simulation results certainly vindicate this expectation.
We would like to note that there are other bootstrap methods that
could have been used, like the wild or residual bootstrap or the
$m$ out of $n$ bootstrap, but it is not clear whether they would yield
consistency at
$\gamma= 1/3$. It would be interesting to explore some of these issues
in future work.
\end{rem}

%s4.1 ###
\subsection{Comparison of two-stage procedures}

By Theorem \ref{Thm_one_stage_IR_estimator_d0},
from the first-stage data,
an asymptotic
$(1-2\beta)$ confidence interval for $d_{0}$ with the true parameter
is given by
\[
\bigl[\hat d_{n_{1}}^{(1)} - Cq_{\beta}n_{1}^{-1/3},
\hat d_{n_{1}}^{(1)} + Cq_{\beta}n_{1}^{-1/3}\bigr]\cap[0,1].
\]
We consider the above
confidence interval as the sampling interval $[L,U]$ with $\gamma=1/3$ and
$K=Cq_{\beta}$.
Then, by Theorem \ref{thm two stage estimator case one LD},
for $f\in\mathscr{F}$ and $\gamma=1/3$,
\[
n^{1/2}\bigl(\tilde
d_{n}^{(2)}-d_{0}\bigr)\stackrel{d}{\rightarrow}C_{2}Z_{1}+C_{3}\mathbb{Z}Z_{2}.
\]
Hence, the corresponding asymptotic $(1-2\alpha)$ confidence
interval of $d_{0}$ is given by
%
%e4.2 ###
%
\begin{equation}\label{ALTSP_2nd_CI}
\bigl[\tilde d_{n}^{(2)} - \tilde q_{\alpha}n^{-1/2}, \tilde d_{n}^{(2)} +
\tilde q_{\alpha}n^{-1/2}\bigr]
\cap[0,1],
\end{equation}
where $\tilde q_{\alpha}$ is the upper
$\alpha$ quantile of $C_{2}Z_{1}+C_{3}\mathbb{Z}Z_{2}$.

%f2 ###
%
\begin{figure}%[b]

\includegraphics{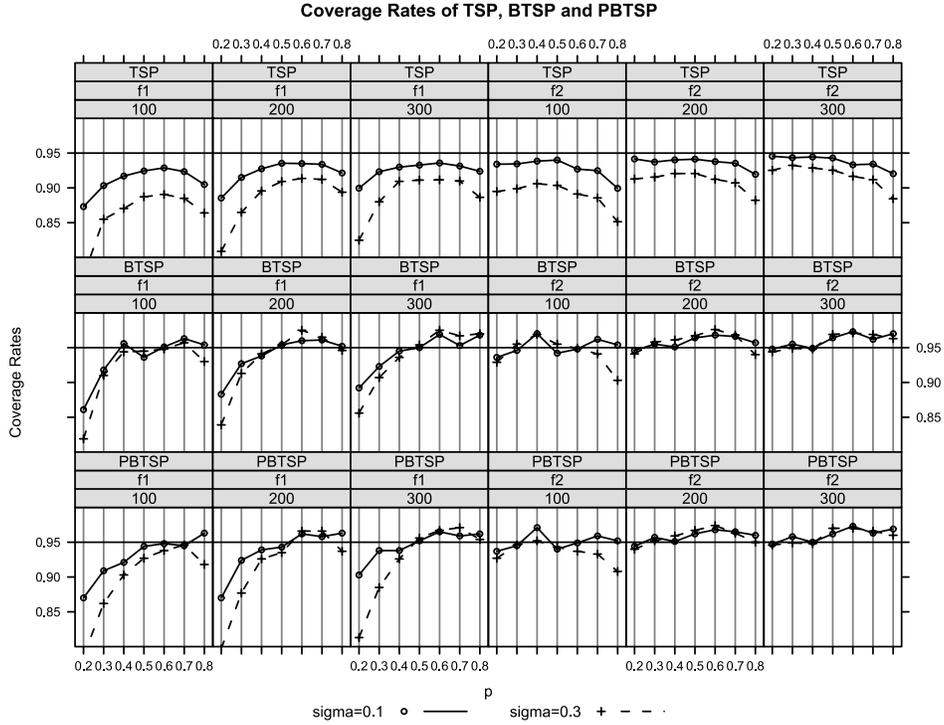}

\caption{Coverage rate plot grouped with different $\sigma$'s.}
\label{Firgure_CR}
\end{figure}

Next, we compare the two-stage procedures, focusing on the coverage rates.
In the first row of Figure \ref{Firgure_CR}, the coverage rates
of the (\ref{ALTSP_2nd_CI}) confidence intervals for combinations of
$f, n$ and $\sigma$ are
shown based on 5000 replications, using the \textit{true} parameters
$f'(d_{0})$ and $\sigma$ (i.e., the true $C$, $C_{2}$ and $C_{3}$
in constructing the confidence intervals). It can be seen that in
general, coverage rates are below the
nominal level 0.95, which is depicted by a solid horizontal line in
each subplot.
This reflects that $\tilde d_{n}^{(2)}$ usually has slow speed of
convergence in distribution.
As expected, the
results improve for small noise levels, larger sample sizes and
functions closer to linearity in
the vicinity of~$d_0$.\looseness=-1

The second row in Figure \ref{Firgure_CR} shows the coverage rates
of the bootstrapped procedure, based on 1000 replicates and 3000
bootstrap samples per replicate,
using the true parameters $f'(d_{0})$ and $\sigma$ at stage one.
It can be seen that the coverage rates achieve the nominal level
with proper first-stage sample proportions,
smaller values of which are preferred
since both average lengths and mean square errors
are usually increasing with $p$ from simulation results not shown in
this paper.
It can be concluded
that the BTSP exhibits superior performance to the TSP
in terms of coverage rates,
especially for settings with $f_{1}$, moderate
noise and relatively small sample sizes.

Finally, the third row in Figure \ref{Firgure_CR} depicts the coverage rates
of the bootstrapped procedure,
when both $f'(d_0)$ and $\sigma$ are estimated from the first-stage data.
The results based on 1000 replicates and 3000 bootstrap samples per
replicate indicate a high level
of agreement with those of the BTSP, which in turn suggests that the PBTSP
is reliable in applications.

Our findings suggest that $p=0.4$ is a good conservative choice for
functions exhibiting
a strong linear trend in the vicinity of $d_0$, while $p=0.5$ is
preferable otherwise.

%s4.2 ###
\subsection{Comparison of one- and two-stage procedures}

We compare next the POSP and the PBTSP, in terms of coverage rates and average
lengths of confidence intervals. We also compare the mean squared
errors of
the first- and second-stage estimates of $d_0$. The results for POSP
are based
on 5000 replications, while those of PBTSP on 1000 replications
and 3000 bootstrap samples per replication, due to its computational
intensity. It can
be seen from Table \ref{Tabel_PABLTSP} that
%
%t1 ###
%
\begin{table}[b]
\caption{CR, AL and MSE stand for coverage rates, average lengths
and mean\break squared errors of PBTSP while \textit{CR1}, \textit{AL1} and \textit{MSE1}
stand for those of POSP.\break ALR and MSER are the ratios of CR over \textit{CR1} and MSE
over \textit{MSE1}, respectively}\label{Tabel_PABLTSP}
\begin{tabular*}{\tablewidth}{@{\extracolsep{\fill}}lccccccccccc@{}}
\hline
$\bolds{f}$ & $\bolds{p}$ & $\bolds{\sigma}$ & $\bolds{n}$ & \textbf{CR} & \textbf{CR1}
& \textbf{AL} & \textbf{AL1} & \textbf{ALR} & \textbf{MSE} &
\textbf{MSE1} & \textbf{MSER} \\
\hline
$f_{1}$ & 0.5 & 0.1 & 100 & \textbf{0.944} & \textbf{0.955} & 0.06 & 0.13 & \textbf
{0.43} & 2e--04 & 1e--03 & \textbf{0.21} \\
& & & 200 & \textbf{0.943} & \textbf{0.953} & 0.04 & 0.10 &
\textbf{0.37} & 1e--04 & 7e--04 & \textbf{0.15} \\
&  & & 300 & \textbf{0.952} & \textbf{0.956} & 0.03 & 0.09
& \textbf{0.35} & 7e--05 & 5e--04 & \textbf{0.14} \\[3pt]
& & 0.3 & 100 & \textbf{0.927} & \textbf{0.942} & 0.21 & 0.27 & \textbf
{0.79} & 3e--03 & 5e--03 & \textbf{0.58} \\
& & & 200 & \textbf{0.935} & \textbf{0.947} & 0.14 & 0.21 &
\textbf{0.63} & 1e--03 & 3e--03 & \textbf{0.39} \\
& & & 300 & \textbf{0.956} & \textbf{0.947} & 0.11 & 0.19 & \textbf
{0.58} & 8e--04 & 2e--03 & \textbf{0.33} \\
[6pt]
$f_{2}$ & 0.4 & 0.1 & 100 & \textbf{0.971} & \textbf{0.966} & 0.06 & 0.16 & \textbf
{0.40} & 2e--04 & 1e--03 & \textbf{0.16} \\
& & & 200 & \textbf{0.951} & \textbf{0.964} & 0.04 & 0.12 &
\textbf{0.34} & 1e--04 & 9e--04 & \textbf{0.13} \\
& & & 300 & \textbf{0.950} & \textbf{0.966} & 0.03 & 0.11
& \textbf{0.31} & 7e--05 & 7e--04 & \textbf{0.11} \\[3pt]
& & 0.3 & 100 & \textbf{0.952} & \textbf{0.948} & 0.24 & 0.32 & \textbf
{0.76} & 5e--03 & 6e--03 & \textbf{0.79} \\
& & & 200 & \textbf{0.959} & \textbf{0.956} & 0.16 & 0.25 &
\textbf{0.62} & 2e--03 & 4e--03 & \textbf{0.46} \\
& & & 300 & \textbf{0.948} & \textbf{0.955} & 0.12 & 0.22 & \textbf
{0.53} & 8e--04 & 3e--03 & \textbf{0.25} \\
\hline
\end{tabular*}
\end{table}
both procedures usually perform well in terms of
coverage rates. Further, under the PBTSP, confidence intervals
usually have shorter average lengths, and the estimates for $d_0$ have
smaller mean squared errors, with slightly more gains accruing in the
$f_2$ case. However, it needs to be pointed out that
both procedures suffer in the case with
large noise and small to moderate sample sizes, especially for $f_{1}$.
\begin{rem}
One of the advantages of the bootstrap procedure,
as~poin\-ted out in Section \ref{Bootstrapping The Two-Stage Estimator},
is that its implementation does not require knowledge of $f'(d_0)$.
One might feel that the practical implementation of the bootstrap
procedure defeats this advantage,
since $f'(d_0)$ is estimated from the first-stage data to construct the
second stage sampling interval.
However, note that only a rough and ready estimate of $f'(d_0)$ would suffice
for the purpose of setting the sampling interval.
On the contrary, to set a confidence interval
directly from the asymptotic distribution of the second-stage estimate
requires a much more precise estimate of $f'(d_0)$.
Thus, the really crucial advantage with the bootstrap is that
it obviates the need for a precise estimate of $f'(d_0)$.
\end{rem}
\begin{rem}
Notice that the sigmoid function $f_{2}$ belongs to class $\mathscr
{F}_2$ for the case
$d_{0}=0.5$, since its second-derivative vanishes at that point. It is
of practical interest to
investigate the performance of the PBTSP for the case where the
regression function at the target point
is close to, but not exactly, linear. We have examined the case for
$f_2$ and $d_0=0.4$
and $0.6$ under the previously considered settings.
The curvatures (i.e., second derivatives) of the regression
functions at these two points are about 0.76 and $-0.76$, respectively.
The results are very close to those obtained for $d_0=0.5$.
\end{rem}
\begin{rem}
In PBTSP, the second stage sampling points $L$ and $U$ are
identified through a Wald-type confidence interval constructed via
estimating $f'(d_0)$ and $\sigma^2$, with $\hat{d}_{n_1}^{(1)}$ at
the center of $[L,U]$. An alternative, albeit ad-hoc way of
obtaining an interval centered at $\hat{d}_{n_1}^{(1)}$ is to set $L
= \hat{d}_{n_1}^{(1)} - L_n/2$ and $U = \hat{d}_{n_1}^{(1)} +
L_n/2$, where $L_n$ is the length of a testing-based confidence
interval for $d_0$ obtained from the first-stage data.
This testing-based interval is obtained as follows: consider testing the
hypothesis $H_{0,d}\dvtx f^{-1}(\theta_0) = d$ vs. $H_{1,d}\dvtx
f^{-1}(\theta_0) \ne d$.
Let $\hat{f}^{(1)}$ denote the usual isotonic estimator of $f$ from
the stage
one data and
$\hat{f}_d^{(1)}$ the constrained isotonic estimator under $H_{0,d}$.
The residual sum of squares based test statistic is given by
\[
\mathrm{RSS}(d) = \frac{\sum_{i=1}^{n_1} (Y_i - \hat{f}_d^{(1)}(x_i))^2 -
\sum_{i=1}^{n_1} (Y_i - \hat{f}^{(1)}(x_i))^2}{\hat{\sigma}^2}  ,
\]
where $\hat{\sigma}^2$ is a consistent estimate of $\sigma^2$. The
inversion procedure assigns $d$ to the confidence set if $\mathrm{RSS}(d)$
falls below an appropriate threshold determined by a
pre-specified quantile of its limit distribution $\mathbb{D}$ (when
$d = d_0$ holds true), which is completely parameter-free and
therefore enables the construction of the confidence set without the
need for nuisance parameter estimation. The limit distribution of
$\mathrm{RSS}(d^0)$ can be derived by
adapting Theorem 2 of \cite{Banerjee2005}
(where a likelihood ratio statistic is dealt with) to the residual sum
of squares statistic in the nonparametric regression setting, but see
also \cite{Banerjee2007d} and \cite{Banerjee2009} for a unified
treatment of likelihood ratio and residual
sum of squares statistics in monotone function problems.

Alternatively, we can use the extremities of the testing-based
confidence interval themselves as the sampling points for
the second stage. For both cases, simulations show that their results
are very similar to
those of PBTSP using the Wald-type confidence interval, thus implying
that the procedure is not
particularly sensitive to the exact specification of $L$ and $U$.
Note that although this testing-based approach has the merit of completely
avoiding the estimation of $f'(d_{0})$, the asymptotic properties of
the corresponding two-stage estimator and its bootstrapped variant
become intractable since neither the testing-based confidence interval
nor the length $L_{n}$ admits an easy analytical characterization,
unlike the analytically simple Wald--type confidence intervals used in
this paper.
To conform to the theoretical development and to save space, we only
present simulation results for such Wald-type stage two sampling intervals.
\end{rem}
\begin{rem}
In the case of $f\in\mathscr{F}_1$, one may question the use of a
linear working model
for approximating $f$ around $d_0$. Instead, fitting a higher order
polynomial working model
may seem more appropriate. We examined the case of $f_1$ using a
quadratic working model.
The results show that this model improves the mean squared error of the
estimates when
the noise is large, but leads to substantial undercoverage.
\end{rem}
\begin{rem}
Our simulation results indicate that good choices for
$p$ are $0.5$ for $f_1$ and 0.4 for $f_2$,
respectively. Our practical recommendation is $p=0.5$,
whenever no prior information about the
linearity of $f$ around $d_0$ is available.
\end{rem}

%s5 ###
\section{Data application}\label{data_application}
We apply our methods to the engineering problem
introduced at the beginning of this paper. We briefly describe the
underlying system next:
consider a complex queueing system comprising $N$
first-in-first-out infinite capacity queues holding different
classes of customers and a set of service resources. These resources
are externally modulated by a stochastic process. The main issue is
to allocate the available resources to the queue in an appropriate
manner so as to maximize the system's throughput. This system
represents a canonical model for wireless data/voice transmissions,
in flexible manufacturing and in call centers (for more details, see
\cite{Bambos2004}).

An important quality of service metric is the average delay of jobs
(over all classes). This quantity can only be obtained through
simulation of the system, due to its analytical intractability. The
average delay of the jobs in a two-class system as a function of its
loading under the optimal throughput policy introduced in
\cite{Bambos2004} is shown in Figure \ref{Figure Data Analysis
Data}. It can be seen that delay is, in general, an increasing
function of the loading. The response was obtained by a discrete
event simulation of the system for each loading, based on 2000
events. Notice that our ability to simulate the system at any
loading in order to obtain the response, allows us to easily
implement the proposed two-stage procedure.

It is of interest to estimate $d_{0}=f^{-1}(\theta_{0})$ for
$\theta_{0}=10$ and 15 units of delay, since around loadings
corresponding to those levels the quality of service provided by the
system exhibits a significant deterioration. For comparing the one- and
two-stage procedures, we fix a budget of $n=82$. A fixed design wth
spacing $0.01$ was used in the interval $[0.14, 0.95]$ to obtain the
one-stage data shown in Figure~\ref{Figure Data Analysis Data} (also in
the left-panel plots of Figure \ref{Figure Data Analysis Compare PBTSP
POSP}). It can be seen that the response is heteroskedastic, but this
does not affect the isotonic regression based estimation of $f$ and
thus of $d_0$. However, it impacts the construction of confidence
intervals through the estimation of the variance at $d_0$. To overcome
this issue, the variance function is estimated locally by the method
proposed in \cite{Muller1987}. More specifically, we compute the
initial local variance estimates with the weights
$(1/\sqrt{2},-1/\sqrt{2})$ and the smoothed variance function by using
\textit{glkerns} in the R package \textit{lokern} with an adaptive
bandwidth, shown in the left panel of Figure \ref{Figure Data Analysis
Estimation Variance}.

%f3 ###
%
\begin{figure}

\includegraphics{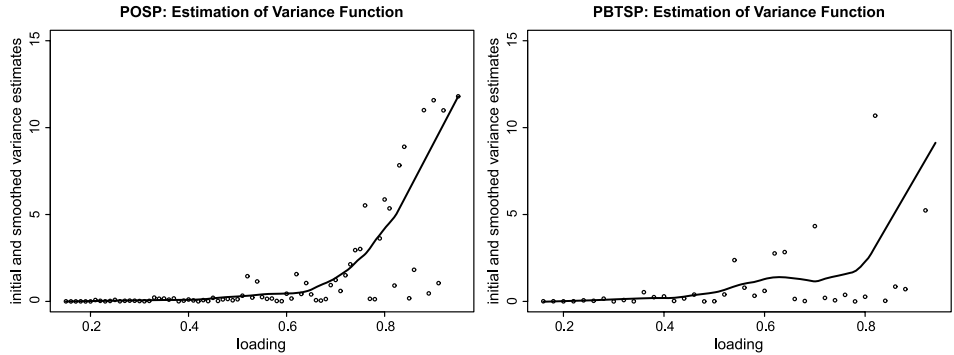}

\caption{Estimation of the variance function in POSP and PBTSP.}
\label{Figure Data Analysis Estimation Variance}
\end{figure}

When implementing the two-stage procedure, we selected every other
point from those used in the one-stage
procedure ($p=0.5$), thus resulting in a fixed design with spacing
$0.02$ on the interval $[0.14,0.94]$.
The initial local variance estimates and smoothed variance function
with the first-stage data
are shown in the right panel of Figure \ref{Figure Data Analysis
Estimation Variance}.
After obtaining the $40=2\times20$ second-stage responses,
the second-stage estimator of $d_{0}$ was computed using
weighted least squares, with weights being the reciprocals of the
estimated local
variances at the corresponding sampling points.

%t2 ###
%
\begin{table}
\caption{Comparing POSP and PBTSP}
\label{Table Data Analysis Compare PBTSP POSP}
\begin{tabular*}{\tablewidth}{@{\extracolsep{\fill}}lccc@{}}
\hline
& & \textbf{POSP} $\bolds{n=82}$ & \textbf{PBTSP} $\bolds{n=81=41+2\times20}$ \\
\hline
$\theta= 10$& estimates of $d_{0}$ & $\hat d_{n}^{(1)}=0.803$ &
$\tilde d_{n}^{(2)}=0.799$ \\
& 95\% CI & $[0.764, 0.841]$ & $[0.794, 0.804]$ \\
[4pt]
$\theta= 15$& estimates of $d_{0}$ & $\hat d_{n}^{(1)}=0.863$ &
$\tilde d_{n}^{(2)}=0.857$ \\
& 95\% CI & $[0.839, 0.887]$ & $[0.845, 0.875]$ \\
\hline
\end{tabular*}
\end{table}

%f4 ###
%
\begin{figure}

\includegraphics{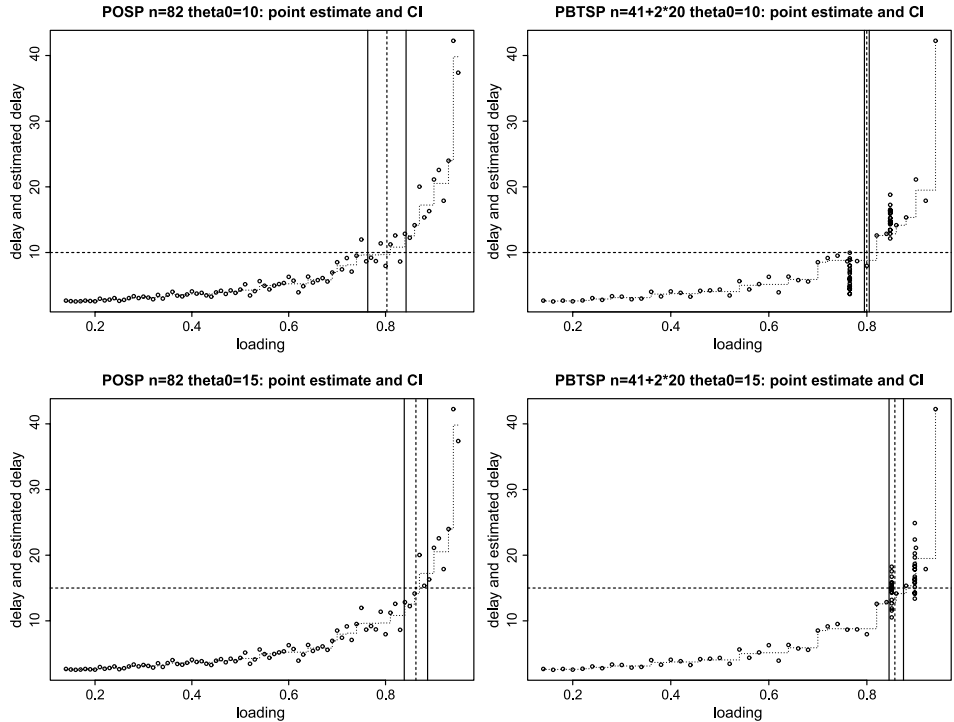}

\caption{Comparing POSP and PBTSP.}
\label{Figure Data Analysis Compare PBTSP POSP}
\end{figure}

The point estimates and the associated 95\% confidence intervals
from the POSP and the PBTSP are given in
Table \ref{Table Data Analysis Compare PBTSP POSP} and
plotted in Figure \ref{Figure Data Analysis Compare PBTSP POSP}.
It can be seen that the point estimates are fairly similar.
More significantly, the confidence intervals from PBTSP are
much shorter than those from POSP,
especially for the case \mbox{$\theta_{0}=10$}.
This can be attributed to two factors: (i) the applicability of the
linear model locally
and (ii) the presence of a strong signal (small noise) for the design
points around~$0.8$.

%%%\input{Conclusions}
%s6 ###
\section{Conclusions}\label{Sec_conclusion}

In this study, a two-stage hybrid procedure for estimating an
inverse regression function at a given point was introduced. The
proposed procedure, by first obtaining a nonparametric estimate of
the regression function and subsequently fitting a parametric linear
model in an appropriately shrinking neighborhood of the parameter of
interest, achieves a $\sqrt{n}$ rate of convergence for the
corresponding estimator. Note that isotonic regression was used in
the first stage as it works with minimal assumptions on the
underlying monotone regression function; nevertheless, other
nonparametric procedures could be used. Further, the local
approximation was primarily based on a linear model, although
quadratic and suitable higher-order approximations could be used,
especially in the presence of a small budget of design points, since
the first stage sampling interval may not be short enough.

A bootstrapped version of the two-stage procedure is provided to
overcome the difficulties posed by the requirement of estimating the
derivative of the regression function at the unknown target point
and the slow speed of convergence, especially with moderate sample
sizes. Its asymptotic properties are also investigated and its
strong consistency established (on this point, see also Remark
\ref{remark:Makowski results}).

Our simulation results indicate that the practical bootstrapped
procedure performs well in a variety of settings. Note that all the
plans can be equipped with random designs for generating the
first-stage data and similar asymptotic results follow.
Nevertheless, for relatively small budgets, fixed designs (e.g.,
quantile based) usually yield improved performance.

Finally, we note that the main results generalize readily to
heteroskedastic models of the form $Y=f(x)+\sigma(x)\varepsilon$, where
$\sigma(x)$ is a scaling function that determines
the error variance. Further, the proposed
procedure should also work for discrete response models; for example,
univariate binary and Poisson regression models with a monotone mean
function. Qualitatively, the results are expected to be analogous to
those established in this study; namely, a $\sqrt{n}$ rate of
convergence would be obtained for the estimator of the parameter of
interest. However, the asymptotic behavior of the second-stage
estimator and its bootstrap counterpart would be different and
depend in an explicit manner on the specific model under
consideration.

\begin{appendix}\label{app}
%s7 ###
\section*{Appendix}
In order to establish the strong consistency of the bootstrapped
two-stage estimator, we need a rate of the almost sure convergence
for the one-stage isotonic
regression estimator $\hat d_{n}^{(1)}$ of $d_0$.
The following\vspace*{1pt} lemma, which is the fixed design version of
Lemma 1 in Durot (2008) \cite{Durot2008}, provides a useful
tail probability for $\hat d_{n}^{(1)}$.
\begin{lem}\label{Lemma:tail probability:fixed design}
Suppose $\mathbb{E}|\varepsilon|^{q}<\infty$ for some $q\geq2$
and \textup{(A2)} and \textup{(A3)} hold.
Then, there exists $K>1$, depending on $q$,
such that for every $\theta\in\mathbb{R}$ and $\eta>0$,
\[
P\bigl(\bigl|\hat d_{n}^{(1)} - d_{0}\bigr|\geq\eta\bigr)\leq K(n\eta^{3})^{-q/2}.
\]
\end{lem}
\begin{pf}
It will be shown that (A2) implies
%
%e7.1 ###
%
\begin{equation}\label{expression inverse distance}
{\sup_{u\in[0,1]}}|F_{n}^{-1}(u)-G^{-1}(x)|\lesssim n^{-1/2}.
\end{equation}
Recall that ``$\lesssim$'' denotes that the left-hand side is less than
a constant times the right-hand side.
Then, reworking the proof of Lemma 1 in Durot \cite{Durot2008} for our
fixed design setting and
an increasing function, and replacing expression (13) in that lemma
with (\ref{expression
inverse distance}) ensures that all subsequent steps go through
yielding the desired conclusion.
To show (\ref{expression inverse distance}), note that
from (A2), we get $|G^{-1}(u)-G^{-1}(v)|\lesssim
|u-v|$ for every $u,v\in[0,1]$. Then
\begin{eqnarray*}
&& {\sup_{u\in[0,1]}}|F_{n}^{-1}(u)-G^{-1}(u)| \\
&&\qquad= \max\{|G^{-1}(G(x_{i}))-G^{-1}(i/n)|,
|G^{-1}(G(x_{i+1}))-G^{-1}(i/n)|, \\
&&\qquad\hspace*{37.7pt} \mbox{for } i=1,2,\ldots, n-1, |G^{-1}(G(x_{1}))|,
|G^{-1}(G(x_{n}))-1|\} \\
&&\qquad\lesssim \max\{
|G(x_{i})-i/n|, |G(x_{i+1})-i/n|\mbox{, for } i=1,2,\ldots, n-1,\\
&&\qquad\hspace*{168.7pt} |G(x_{1})-0|,
|G(x_{n})-1|\} \\
&&\qquad= {\sup_{x\in[0,1]}}|F_{n}(x)-G(x)|
\end{eqnarray*}
gives (\ref{expression inverse distance}) again by (A2).
\end{pf}

With the help of Lemma \ref{Lemma:tail probability:fixed design},
next we show that $n^{1/3}$ is a boundary rate of almost sure convergence.
\begin{lem}\label{Lemma LIL one stage estimator}
If \textup{(A2)} to \textup{(A4)} hold, for each $\alpha>0$,
\[
P\Bigl(\lim_{n\rightarrow\infty} n^{1/3-\alpha}
\bigl|\hat d_{n}^{(1)}-d_{0}\bigr| = 0 \Bigr)= 1.
\]
Thus, for every $r<1/3$,
$\lim_{n\rightarrow\infty}n^{r}(\hat d_{n}^{(1)}-d_{0}) =0, (P\mbox
{-a.s.})$.
\end{lem}
\begin{pf}
Use the notations $K$, $q$ and $\eta$ in Lemma \ref{Lemma:tail
probability:fixed design}. Denote $K'=K\eta^{-3q/2}$ and
$A_{n}=\{n^{1/3-\alpha}|\hat d_{n}^{(1)}-d_{0}|\geq\eta\}$.
By Lemma \ref{Lemma:tail probability:fixed design},
$
P(A_{n})\leq\break K'n^{-3\alpha q/2}
$ for each $\alpha>0$. On the other hand, (A4) allows $q$ to be
arbitrarily large. Choosing $q>2/(3\alpha)$ gives $
\sum_{n=1}^{\infty}P(A_{n})
\leq K'\sum_{n=1}^{\infty} n^{-3\alpha q/2}
<\infty
$.
Note that $\eta>0$ is arbitrary.
Therefore,
$n^{1/3-\alpha}
|\hat d_{n}^{(1)}-d_{0}|$
converges to $0$ almost surely
(see Corollary on pages 254 and 255 in Shiryaev \cite{Shiryaev1995}),
which completes the proof.
\end{pf}
\begin{rem}
Note that
Lemmas \ref{Lemma:tail probability:fixed design}
and \ref{Lemma LIL one stage estimator} hold for not only sequences,
but also triangular arrays of design points and random errors.
\end{rem}
\begin{rem}\label{remark:Makowski results}
The proof of Lemma \ref{Lemma LIL one stage estimator} implies
$n^{1/3-\alpha}(\hat d_{n}^{(1)}-d_{0})\stackrel{\mathrm{a.s.}}{\rightarrow} 0$
for each $\alpha\in(0,1/3)$ given $q>2/(3\alpha)$.
Then, $\mathbb{E}|\varepsilon|^{8}<\infty$ ensures
$n^{\beta}(\hat d_{n}^{(1)}-d_{0})\stackrel{\mathrm{a.s.}}{\rightarrow} 0$ for
each $\beta<1/4$. However, this almost sure convergence result
actually holds under a weaker condition
$\mathbb{E}|\varepsilon|^{3}<\infty$ by theorem in Makowski
\cite{Makowski1975} and Remark 4 in Makowski \cite{Makowski1973}.
This shows that it might be possible to weaken the assumption (A4) a
little. Essentially, it means that it might be possible to weaken
the condition on the random error in Lemma \ref{Lemma:tail
probability:fixed design}. In fact, this possibility has been
mentioned in
Durot's papers on isotonic regression
\cite{Durot2002,Durot2007,Durot2008}.
\end{rem}

%%%\input{Proofs_TwoStagePlan}
%s7.1 ###
\subsection{\texorpdfstring{Proofs for results in Section \protect\ref{Two-Stage Estimator}}%
{Proofs for results in Section 3.1}}
For the simplicity of notation, from now on
denote $\delta_{d}=\hat d_{n_{1}}^{(1)}-d_{0}$,
$\varepsilon_{i}^{+}=\varepsilon''_{i}+\varepsilon'_{i}$,
$\varepsilon_{i}^{-}=\varepsilon''_{i}-\varepsilon'_{i}$,
$f_{UL}^{+}=f(U)+f(L)$, $f_{UL}^{-}=f(U)-f(L)$,
$R_{UL}^{+}=R_{U}+R_{L}$, $R_{UL}^{-}=R_{U}-R_{L}$,
$R_{UL}^{\prime+}=R_{U}'+R_{L}'$, $R_{UL}^{\prime-}=R_{U}'-R_{L}'$.
Recall $Y_{i}^{+}=Y''_{i}+Y'_{i}$ and
$Y_{i}^{-}=Y''_{i}-Y'_{i}$.
\begin{pf*}{Proof of Lemma
\ref{Two_Ends_General_Case_One_Lemma_Consistency}}
Consider the following Taylor expansions:
%
%e7.2 ###
%
\begin{eqnarray}\label{Two_Ends_Talor_3_prime_1}
f(U) &=& f\bigl(\hat d_{n_{1}}^{(1)}+Kn_{1}^{-\gamma}\bigr)\nonumber\\[-2pt]
&=& f(d_{0})+f'(d_{0})(\delta_{d}+Kn_{1}^{-\gamma})\\[-2pt]
&&{}+ (1/2)f''(d_{0})(\delta_{d}+Kn_{1}^{-\gamma})^{2} + R_{U},\nonumber[-2pt]
\\
%
%e7.3 ###
%
\label{Two_Ends_Talor_3_prime_2}
f(L) &=& f\bigl(\hat d_{n_{1}}^{(1)}-Kn_{1}^{-\gamma}\bigr)\nonumber\\[-2pt]
&=& f(d_{0})+f'(d_{0})(\delta_{d}-Kn_{1}^{-\gamma})\\[-2pt]
&&{}+ (1/2)f''(d_{0})(\delta_{d}-Kn_{1}^{-\gamma})^{2}+ R_{L},\nonumber
\end{eqnarray}
where
$
R_{U}=f'''(\xi_{1})(\delta_{d}+Kn_{1}^{-\gamma})^{3}/6$, $ R_{L}=
f'''(\xi_{2})(\delta_{d}-Kn_{1}^{-\gamma})^{3}/6$,
$\xi_{1}$ lies between $d_{0}$ and $\hat
d_{n_{1}}^{(1)}+Kn_{1}^{-\gamma}$
and $\xi_{2}$ lies between $d_{0}$ and $\hat
d_{n_{1}}^{(1)}-Kn_{1}^{-\gamma}$.
Since $\hat d_{n_{1}}^{(1)}$ converges to $d_{0}$ in probability by
Theorem \ref{Thm_one_stage_IR_estimator_d0},
so do $\xi_{1}$ and $\xi_{2}$.

Then, from (\ref{TED G beta0 beta1}), the definitions of $Y'_{i}$
and $Y''_{i}$ and the Taylor expansions
(\ref{Two_Ends_Talor_3_prime_1}) and
(\ref{Two_Ends_Talor_3_prime_2}), we get\vspace{-4pt}
\begin{eqnarray*}
\hat\beta_{1}
&=& (2Kn_{1}^{-\gamma}n_{2})^{-1}\sum_{i=1}^{n_{2}}Y_{i}^{-}
= (2Kn_{1}^{-\gamma})^{-1}f_{UL}^{-}
+ (2Kn_{1}^{-\gamma}n_{2})^{-1}\sum_{i=1}^{n_{2}}\varepsilon_{i}^{-}\\[-4pt]
&=& f'(d_{0}) + f''(d_{0})\delta_{d}
+ (2Kn_{1}^{-\gamma})^{-1}R_{UL}^{-}
+ (2Kn_{1}^{-\gamma}n_{2})^{-1}\sum_{i=1}^{n_{2}}\varepsilon_{i}^{-}.
\vspace{-5pt}
\end{eqnarray*}

From Theorem \ref{Thm_one_stage_IR_estimator_d0}, $
\delta_{d}\stackrel{P}{\rightarrow} 0
$; and by the Lindeberg--Feller CLT for triangular arrays,
for $\gamma
\in(0,1/2)$, $
(n_{1}^{\gamma}/n_{2})\sum_{i=1}^{n_{2}}\varepsilon_{i}^{-}
\stackrel{P}{\rightarrow}0
$. Next we show that
$
R_{UL}^{-}/(2Kn_{1}^{-\gamma})
\stackrel{P}{\rightarrow} 0
$ for $\gamma\in(0,1)$.
Hence, for $\gamma\in(0,1/2)$ we get
$
\hat\beta_{1} \stackrel{P}{\rightarrow} f'(d_{0})
$.
It suffices to show both
$n_{1}^{\gamma}R_{U}$ and $n_{1}^{\gamma}R_{L}$ converge to 0 in
probability for $\gamma\in(0,1)$.
We only show the former; the latter follows in an analogous
manner.\looseness=-1

From the definition of $R_U$, we have
%
%e7.4 ###
%
\begin{eqnarray}\label{Result_n1gamma_ru_expansion}
n_{1}^{\gamma}R_{U}
\nonumber&=& (1/6)n_{1}^{\gamma}f'''(\xi_{1})(\delta
_{d}+Kn_{1}^{-\gamma})^{3}\nonumber\\[-8pt]\\[-8pt]
&=& (1/6)f'''(\xi_{1})
[n_{1}^{\gamma} \delta_{d}^{3}
+ 3K\delta_{d}^{2}
+ 3K^{2}n_{1}^{-\gamma}\delta_{d}
+ K^{3}n_{1}^{-2\gamma}].\nonumber
\end{eqnarray}
Theorem \ref{Thm_one_stage_IR_estimator_d0} coupled with Slutsky's
lemma, shows that the sum of the four terms within the square
bracket in (\ref{Result_n1gamma_ru_expansion}) is $o_{P}(1)$ for
$\gamma\in(0,1)$. Thus, we have
$n_{1}^{\gamma}R_{U}=f'''(\xi_{1})o_{P}(1)$. Since
$f'''(\cdot)$ is uniformly bounded around $d_{0}$ and
$\xi_{1}\rightarrow d_{0}$ in probability,
$f'''(\xi_{1})o_{P}(1)=o_{P}(1)$. This shows that
$n_{1}^{\gamma}R_{U}$ converges to~$0$ in probability
for $\gamma\in(0,1)$. Obviously, $R_{U}=o_{P}(1)$.

Then, for $\gamma\in(0,1/2)$,
\begin{eqnarray*}
\hat\beta_{0}
&=& (2n_{2})^{-1}\sum_{i=1}^{n_{2}}Y_{i}^{+}-\hat d_{n_{1}}^{(1)} \hat
\beta_{1}\\
&=& f(d_{0})
+ (1/2)f''(d_{0})[ \delta_{d}^{2} + K^{2}n_{1}^{-2\gamma} ]\\
&&{}+ f'(d_{0})\delta_{d} + (1/2)R_{UL}^{+}
+ (2n_{2})^{-1}\sum_{i=1}^{n_{2}}\varepsilon_{i}^{+}
- \hat d_{n_{1}}^{(1)} \hat
\beta_{1}\\
&\stackrel{P}{\rightarrow}& f(d_{0}) - d_{0}f'(d_{0}).
\end{eqnarray*}

Finally, for $\gamma\in(0,1/2)$, the weak consistency of $\hat
\beta_{1}$ and $\hat\beta_{0}$ gives
$
\tilde d_{n}^{(2)}
= (\theta_{0}-\hat\beta_{0})/(\hat\beta_{1})
\stackrel{P}{\rightarrow} d_{0}
$.
\end{pf*}
\begin{pf*}{Proof of Theorem \ref{thm two stage estimator case one LD}}
First, suppose $f\in\mathscr{F}_{1}$.
From (\ref{TED_G_tilded}), the definitions of $Y'_{i}$ and $Y''_{i}$
and the Taylor expansions (\ref{Two_Ends_Talor_3_prime_1}) and
(\ref{Two_Ends_Talor_3_prime_2}), we get
\begin{eqnarray*}
\tilde d_{n}^{(2)} - d_{0}
&=& (1/\hat\beta_{1})\Biggl[f(d_{0})-(2n_{2})^{-1}\sum_{i=1}^{n_{2}}Y_{i}^{+}\Biggr]
+ \delta_{d} \\
&=& \bigl(1/f'(d_{0})\bigr)
\Biggl[ f(d_{0})- (2n_{2})^{-1}\sum_{i=1}^{n_{2}}Y_{i}^{+} \Biggr] +
\delta_{d}\\
&&{} + (f'(d_{0})\hat\beta_{1})^{-1}\bigl(f'(d_{0})-\hat\beta_{1}\bigr)
\Biggl[ f(d_{0})- (2n_{2})^{-1}\sum_{i=1}^{n_{2}}Y_{i}^{+} \Biggr]\\
&=& S_{1}+S_{2}\times S_{3},
\end{eqnarray*}
where
\begin{eqnarray*}
S_{1} &=& -f''(d_{0})(2f'(d_{0}))^{-1}(\delta_{d}^{2}
+ K^{2}n_{1}^{-2\gamma} )\\
& &{} - (2f'(d_{0}))^{-1}R_{UL}^{+}
- (2f'(d_{0})n_{2})^{-1}\sum_{i=1}^{n_{2}}\varepsilon_{i}^{+},\\
S_{2} &=& (f'(d_{0})\hat\beta_{1})^{-1}
\Biggl[ f''(d_{0})\delta_{d}+ (2Kn_{1}^{-\gamma})^{-1}R_{UL}^{-}
+(2Kn_{1}^{-\gamma}n_{2})^{-1}\sum_{i=1}^{n_{2}}\varepsilon_{i}^{-}
\Biggr],\\
S_{3} &=& f'(d_{0})\delta_{d}
+ (1/2)f''(d_{0})(\delta_{d}^{2} + K^{2}n_{1}^{-2\gamma} )
+ (1/2)R_{UL}^{+}
+ (2n_{2})^{-1}\sum_{i=1}^{n_{2}}\varepsilon_{i}^{+}.
\end{eqnarray*}

Next, consider the exact stochastic orders of the terms $S_1,
S_2$ and $S_3$. We start with $S_1$. From Theorem
\ref{Thm_one_stage_IR_estimator_d0},
$\delta_{d}^{2}=O_{P}(n^{-2/3})$; for $\gamma>0$,
$n_{1}^{-2\gamma}=O_{P}(n^{-2\gamma})$,
$R_{U}=O_{P}(n^{-1})+O_{P}(n^{-3\gamma})$,
$R_{L}=O_{P}(n^{-1})+O_{P}(n^{-3\gamma})$ and\break
$n_{2}^{-1}\sum_{i=1}^{n_{2}}\varepsilon_{i}^{+}=O_{P}(n^{-1/2})$.
Note that these are the exact rates of weak convergence. Then, for
$\gamma\in(0,1/2)$,
$
S_{1}= T_{1}+T_{2}
+ o_{P}(n^{-2\gamma}\vee n^{-1/2}),
$
where
\[
T_{1}
= -(2f'(d_{0}))^{-1}f''(d_{0})K^{2}n_{1}^{-2\gamma},\qquad
T_{2} = - (2f'(d_{0})n_{2})^{-1}\sum_{i=1}^{n_{2}}
\varepsilon_{i}^{+}.
\]
Thus, the possible main terms of $S_{1}$ are $T_{1}$ and $T_{2}$.
In the same way, we can obtain the main terms of $S_{2}$ and $S_{3}$
and then those of $S_{2}\times S_{3}$. Finally, we have
$
S_{1}+S_{2}\times S_{3}
= T_{1}+T_{2}+T_{3}+R
$,
where
\[
T_{3} = (2K\hat\beta_{1}n_{1}^{-\gamma}n_{2})^{-1}
\delta_{d}\sum_{i=1}^{n_{2}}\varepsilon_{i}^{-},\qquad
R = o_{P}(n^{-2\gamma}\vee n^{-1/2}\vee n^{\gamma-5/6}).
\]
It is easy to see that among the three rates
$n^{-2\gamma}$, $n^{-1/2}$ and $n^{\gamma-5/6}$,
the first, second or last one is slowest
according as $\gamma$ belongs to the interval $(1,1/4)$,
$(1/4,1/3)$, or $(1/3,1/2)$, respectively; the first and the second are
the slowest
for $\gamma=1/4$; while the second and the last ones are the slowest
for $\gamma=1/3$.
In other words, $T_{1}$, $T_{2}$ or $T_{3}$ becomes the main term
according as $\gamma\in(0,1/4)$, $\gamma\in(1/4,1/3)$ or $\gamma
\in
(1/3,1/2)$, respectively. When $\gamma=1/4$, both $T_{1}$ and $T_{2}$
become the
main terms and when $\gamma=1/3$, both $T_{2}$ and $T_{3}$ become
the main terms.

Then, by Theorem \ref{Thm_one_stage_IR_estimator_d0}, the
Lindeberg--Feller CLT
for triangular arrays, Slutsky's lemma and the Continuous
Mapping theorem, and noting that $n_{1}^{1/3}\delta_{d}$ is
independent of $n_{2}^{-1/2}\sum_{i=1}^{n_{2}}\varepsilon_{i}^{+}$ and
$n_{2}^{-1/2}\sum_{i=1}^{n_{2}}\varepsilon_{i}^{-}$ and that
$\varepsilon_{i}^{+}$ is uncorrelated with $\varepsilon_{i}^{-}$, we
obtain the results of the five cases for $f\in\mathscr{F}_{1}$
defined by the different ranges of $\gamma$ in the statement of the
theorem.

For the purpose of illustration, we outline the case
$\gamma=1/3$, for which $T_{2}+T_{3}$ is the main term with exact
stochastic order $O_{P}(n^{-1/2})$. Thus $n^{1/2}(\tilde
d_{n}^{(2)} - d_{0})$ and $n^{1/2}(T_{2}+T_{3})$ have the same
asymptotic distribution. Since
\[
\Biggl( n_{1}^{1/3}\delta_{d},
n_{2}^{-1/2}\sum_{i=1}^{n_{2}}
\varepsilon_{i}^{+},
n_{2}^{-1/2}\sum_{i=1}^{n_{2}}
\varepsilon_{i}^{-}\Biggr)
\stackrel{d}{\rightarrow}
(C\mathbb{Z}, c Z_{1}, c Z_{2}),
\]
where $\mathbb{Z}$ follows Chernoff
distribution, independent of $Z_{1}$, $Z_{2}$ which are i.i.d. $N(0,1)$,
and $c=\sqrt{2}\sigma$,
by Continuous Mapping
theorem, we have
\[
n^{1/2}(T_{2}+T_{3})
\stackrel{d}{\rightarrow}
- C_{2}Z_{1} +
(1/K)C_{2}C\mathbb{Z}Z_{2}.
\]
Note that $-C_{2}Z_{1}$ can be replaced by
$C_{2}Z_{1}$ since $N(0,1)$ and $-N(0,1)$ have the same distribution.
In similar fashion, we obtain the asymptotic results for the other
four cases.

Carefully examining the above proof reveals that
the conclusions with $\gamma\in(1/4,1/2)$
also hold for $f\in\mathscr{F}$.
Thus, it remains to show the cases $f\in\mathscr{F}_{2}$
and $\gamma\in(1/8,1/4]$.

For $f\in\mathscr{F}_{2}$, consider the following Taylor expansions:
%
%e7.5 ###
%
\begin{eqnarray}\label{Two_Ends_Case_Two_Talor_f_U}
f(U) &=& f\bigl(\hat d_{n_{1}}^{(1)}+Kn_{1}^{-\gamma}\bigr)\nonumber\\
&=& f(d_{0})+f'(d_{0})(\delta_{d}+Kn_{1}^{-\gamma})
\\
&&{}+
(1/6)f'''(d_{0})(\delta_{d}+Kn_{1}^{-\gamma})^{3}+R_{U}',\nonumber\\
%
%e7.6 ###
%
\label{Two_Ends_Case_Two_Talor_f_L}
f(L) &=& f\bigl(\hat d_{n_{1}}^{(1)}-Kn_{1}^{-\gamma}\bigr)\nonumber\\
&=&
f(d_{0})+f'(d_{0})(\delta_{d}-Kn_{1}^{-\gamma})\\
&&{} +
(1/6)f'''(d_{0})(\delta_{d}-Kn_{1}^{-\gamma})^{3}+R_{L}',\nonumber
\end{eqnarray}
where
$
R_{U}'=f^{(4)}(\xi_{1})(\delta_{d}+Kn_{1}^{-\gamma})^{4}/24
$,
$
R_{L}'=f^{(4)}(\xi_{2})(\delta_{d}-Kn_{1}^{-\gamma})^{4}/24
$,
$\xi_{1}$ lies between $d_{0}$ and $\hat
d_{n_{1}}^{(1)}+Kn_{1}^{-\gamma}$
and $\xi_{2}$ lies between $d_{0}$ and $\hat
d_{n_{1}}^{(1)}-Kn_{1}^{-\gamma}$.

Then, for
$\gamma\in(1/8,1/2)$,
\begin{eqnarray*}
\tilde d_{n}^{(2)} - d_{0}
&=& (1/\hat\beta_{1})\Biggl[f(d_{0})-(2n_{2})^{-1}\sum_{i=1}^{n_{2}}Y_{i}^{+}\Biggr]
+ \delta_{d} \\
&=& \bigl(1/f'(d_{0})\bigr)
\Biggl[ f(d_{0})- (2n_{2})^{-1}\sum_{i=1}^{n_{2}}Y_{i}^{+} \Biggr]+
\delta_{d}\\
&&{} + (f'(d_{0})\hat\beta_{1})^{-1}\bigl(f'(d_{0})-\hat\beta_{1}\bigr)
\Biggl[ f(d_{0})- (2n_{2})^{-1}\sum_{i=1}^{n_{2}}Y_{i}^{+} \Biggr]\\
&=& S_{1} + S_{2}\times S_{3},
\end{eqnarray*}
where
\begin{eqnarray*}
S_{1} &=& -(6f'(d_{0}))^{-1}f'''(d_{0})\delta_{d}^{3}
- (2f'(d_{0}))^{-1}f'''(d_{0})\delta_{d} K^{2}n_{1}^{-2\gamma}\\
& &{} - (2f'(d_{0}))^{-1}R_{UL}^{\prime+}
- (2f'(d_{0})n_{2})^{-1}\sum_{i=1}^{n_{2}}\varepsilon_{i}^{+},\\
S_{2} &=& (f'(d_{0})\hat\beta_{1})^{-1}
\Biggl[(1/2)f'''(d_{0})\delta_{d}^{2}
+(1/6)f'''(d_{0})K^{2}n_{1}^{-2\gamma}\\
& &\hspace*{64.4pt}{} + (2Kn_{1}^{-\gamma})^{-1}R_{UL}^{\prime-}
+ (2Kn_{1}^{-\gamma}n_{2})^{-1}\sum_{i=1}^{n_{2}}\varepsilon
_{i}^{-}\Biggr],\\
S_{3} &=& \Biggl\{ f'(d_{0})\delta_{d}
+ \frac{f'''(d_{0})}{6}\delta_{d}^{3}
+ \frac{f'''(d_{0})}{2}\delta_{d} K^{2}n_{1}^{-2\gamma}
+\frac{1}{2}R_{UL}^{\prime+}
+ \frac{1}{2n_{2}}\sum_{i=1}^{n_{2}}\varepsilon_{i}^{+}\Biggr\}.
\end{eqnarray*}

Similar to the previous argument on the exact weak convergence rates,
$S_{1}+S_{2}\times S_{3}=T_{1}+T_{2}+R'$, where
\[
T_{1}
= - (2f'(d_{0})n_{2})^{-1}\sum_{i=1}^{n_{2}}\varepsilon_{i}^{+},\qquad
T_{2}
= (1/\hat\beta_{1})\delta_{d}
(2Kn_{1}^{-\gamma}n_{2})^{-1}\sum_{i=1}^{n_{2}}\varepsilon_{i}^{-},
\]
and $R'$ is the sum of the remaining terms which converges to $0$
faster than $T_{1}$ and~$T_{2}$. Then $T_{1}$ becomes the main term
for $\gamma\in(1/8, 1/3)$ and the weak convergence result for
$f\in\mathscr{F}_{2}$ and $\gamma\in(1/8,1/4]$ follows easily from
the Lindeberg--Feller central limit theorem for triangular arrays
and Slutsky's lemma. This completes the proof.
\end{pf*}

%%%\input{Proofs_BootstrappedTwoStagePlan}
%s7.2 ###
\subsection{\texorpdfstring{Proofs for results in Section \protect\ref{Bootstrapped Two-Stage Estimator}}%
{Proofs for results in Section 3.2}}
To simplify arguments, we introduce a notation on the rate of almost
sure convergence.
Suppose $\{\zeta_{n}\}$ is a sequence of random variables and
$b\in\mathbb{R}$. Write $\zeta_{n}=B_{as}(b)$ if
$n^{\alpha}\zeta_{n}$ converges to 0 almost surely for every
$\alpha<b$. It is easy to verify that
$B_{as}(b_{1})+B_{as}(b_{2})=B_{as}(b_{1})$ and
$B_{as}(b_{1})B_{as}(b_{2})=B_{as}(b_{1}+b_{2})$ if $b_{1}\leq
b_{2}\in\mathbb{R}$.
Note that $\zeta_{n}=B_{as}(b)$ for some $b>0$ implies $\zeta
_{n}\rightarrow0$ almost surely.
Denote
$V_{i}^{+}\equiv\varepsilon_{i}^{\star+}=\varepsilon_{i}^{\prime\prime\star
}+\varepsilon_{i}^{\prime\star}$ and
$V_{i}^{-}\equiv\varepsilon_{i}^{\star-}=\varepsilon_{i}^{\prime\prime\star
}-\varepsilon_{i}^{\prime\star}$.
Recall
$Y_{i}^{\star+}=Y_{i}^{\prime\prime\star}+Y_{i}^{\prime\star}$ and
$Y_{i}^{\star-}=Y_{i}^{\prime\prime\star}-Y_{i}^{\prime\star}$.
\begin{pf*}{Proof of Lemma \ref{Lemma_Bootstrapped_Root_Consistency}}
The proof of Lemma
\ref{Two_Ends_General_Case_One_Lemma_Consistency} establishes the
weak consistency of $\hat\beta_{1}$ for the case $\gamma\in
(0,1/2)$. In fact, under the setting of the bootstrapped two-stage
procedure, the strong consistency of $\hat\beta_{1}$ can be
obtained.

From the proof of Lemma
\ref{Two_Ends_General_Case_One_Lemma_Consistency}, it suffices to
show $\delta_{d}$,
$(n_{1}^{\gamma}/n_{2})\sum_{i=1}^{n_{2}}\varepsilon_{i}^{-}$
and $R_{UL}^{-}/(2Kn_{1}^{-\gamma})$
converge to $0$ almost surely. Lemma
\ref{Lemma LIL one stage estimator} shows that $\delta_{d}$
converges to $0$ almost surely, while Lemma \ref{Lemma_LOIL_Mean}
establishes that
$(n_{1}^{\gamma}/n_{2})\sum_{i=1}^{n_{2}}\varepsilon_{i}^{-}$ converges
to $0$ almost surely for $\gamma\in(0,1/2)$. Thus, it suffices to
show that both $n_{1}^{\gamma}R_{U}$ and $n_{1}^{\gamma}R_{L}$
converge to 0 almost surely for $\gamma\in(0,1)$. Next, we show the
former; the latter follows analogously.

Since $\xi_{1}$ lies between $d_{0}$ and $\hat
d_{n_{1}}^{(1)}+Kn_{1}^{-\gamma}$ and the latter converges to
$d_{0}$ almost surely, we know $\xi_{1}$ converges to $d_{0}$ almost
surely. On the other hand, $f'''(\cdot)$ is uniformly bounded around
$d_{0}$; thus, $f'''(\xi_{1})$ is almost surely bounded.
Further, by Lemma \ref{Lemma LIL one stage estimator},
the four terms within square brackets on the right-hand side of (\ref
{Result_n1gamma_ru_expansion})
are $B_{as}(1-\gamma)$, $B_{as}(2/3)$, $B_{as}(1/3+\gamma)$ and
$B_{as}(2\gamma)$.
Thus,
$n_{1}^{\gamma}R_{U}$ almost surely converges to $0$ for
$\gamma\in(0,1)$.

So, for $\gamma\in(0,1/2)$, we have
$
\hat\beta_{1} \rightarrow f'(d_{0}), (P\mbox{-a.s.})
$.

Next, we establish the conditional weak consistency of $\hat
\beta_{1}^{\star}$ for $f\in\mathscr{F}$.
From (\ref{Bootstrap 2nd stage estimator beta
0 beta1}), we get
\[
\hat\beta_{1}^{\star}
=(2Kn_{1}^{-\gamma}n_{2})^{-1}\sum_{i=1}^{n_{2}}Y_{i}^{\star-}
= T_{1} + T_{2},
\]
where
\[
T_{1} = (2Kn_{1}^{-\gamma}n_{2})^{-1}\sum_{i=1}^{n_{2}}\varepsilon
_{i}^{\star-},\qquad
T_{2} = (2Kn_{1}^{-\gamma})^{-1}f_{UL}^{-}.
\]

Hence, we have
$
T_{1}
= T_{11}+T_{12}
$, where
\[
T_{11} = s(2Kn_{1}^{-\gamma}n_{2})^{-1}
\sum_{i=1}^{n_{2}} (V_{i}^{-}-\nu^{-})/s ,\qquad
T_{12} = (2Kn_{1}^{-\gamma}n_{2})^{-1}\sum_{i=1}^{n_{2}}\varepsilon_{i}^{-},
\]
$
V_{i}^{-}=\varepsilon_{i}^{\star-},
\nu^{-}=E_{\star}[V_{i}^{-}]=(1/n_{2})\sum_{i=1}^{n_{2}}\varepsilon_{i}^{-},
$
and
\[
s^{2}=\operatorname{Var}_{\star}[V_{i}^{-}]=\frac{1}{n_{2}}\sum
_{i=1}^{n_{2}}(\varepsilon''_{i})^{2}
-\Biggl(\frac{1}{n_{2}}\sum_{i=1}^{n_{2}}\varepsilon''_{i}\Biggr)^{2}
+\frac{1}{n_{2}}\sum_{i=1}^{n_{2}}(\varepsilon'_{i})^{2}
-\Biggl(\frac{1}{n_{2}}\sum_{i=1}^{n_{2}}\varepsilon'_{i}\Biggr)^{2}.
\]

For $\gamma\in(0,1/2)$, gives that
\mbox{$
T_{12}\rightarrow0$}, $(P\mbox{-a.s.})
$
by Lemma \ref{Lemma_LOIL_Mean}
and
\mbox{$
T_{11}\stackrel{P^{\star}}{\rightarrow} 0$}, $(P\mbox{-a.s.})
$
by Lemma \ref{Lemma_TCLT_V_Plus_Minus}
and Slutsky's lemma.
Thus, $
T_{1}\stackrel{P^{\star}}{\rightarrow} 0, (P\mbox{-a.s.})
$.

Next, we consider $T_{2}$. By the almost sure convergence
of $\delta_{d}$
and $R_{UL}^{-}/\break (2Kn_{1}^{-\gamma})$,
we have, for $\gamma\in(0,1)$,
\[
T_{2}
= f'(d_{0}) + f''(d_{0})\delta_{d}
+ (2Kn_{1}^{-\gamma})^{-1} R_{UL}^{-}
\rightarrow f'(d_{0}),\qquad  (P\mbox{-a.s.}).
\]
Thus, for $f\in\mathscr{F}$ and $\gamma\in(0,1/2)$,
$
T_{2}\rightarrow f'(d_{0}), (P\mbox{-a.s.})
$. Therefore, we get $\hat\beta_{1}^{\star}
\stackrel{P^{\star}}{\rightarrow} f'(d_{0}), (P\mbox{-a.s.})$.
\end{pf*}
\begin{pf*}{Proof of Theorem \ref{Thm Bootstrapped two stage
estimator LD}}
From (\ref{TED_G_tilded}) and (\ref{Bootstrap_2nd_stage_estimator_d0}),
\[
n^{1/2}\bigl(\tilde d_{n}^{(2)\star}-\tilde d_{n}^{(2)}\bigr)
= -T_{1} + T_{2},
\]
where
\begin{eqnarray*}
T_{1}&=&(f'(d_{0})2n_{2})^{-1}n^{1/2}
\sum_{i=1}^{n_{2}}
( Y_{i}^{\star+} -
Y_{i}^{+} ),\\
T_{2} &=&n^{1/2} \Biggl[ \bigl(1/\hat\beta_{1}^{\star
}-1/f'(d_{0})\bigr)
\Biggl( f(d_{0}) - (2n_{2})^{-1}\sum_{i=1}^{n_{2}}Y_{i}^{\star+}
\Biggr)\\
&&\hspace*{23.8pt}{} -\bigl(1/\hat\beta_{1}-1/f'(d_{0})\bigr)
\Biggl( f(d_{0}) - (2n_{2})^{-1}\sum_{i=1}^{n_{2}}Y_{i}^{+} \Biggr)\Biggr].
\end{eqnarray*}

By the definitions of $Y'_{i}, Y''_{i}, Y'^{\star}_{i},
Y''^{\star}_{i}$,
\[
T_{1}
= n^{1/2}(f'(d_{0})2n_{2})^{-1}
\sum_{i=1}^{n_{2}}
(\varepsilon_{i}^{\star+} -
\varepsilon_{i}^{+} )
= sn^{1/2}(2f'(d_{0})n^{1/2})^{-1}
\sum_{i=1}^{n_{2}}\frac{V_{i}^{+}-\nu^{+}}{s\sqrt{n_{2}}},
\]
where
\[
V_{i}^{+}=\varepsilon_{i}^{\star+},\qquad
\nu^{+}=E_{\star}[V_{i}^{+}
]=(1/n_{2})\sum_{i=1}^{n_{2}}\varepsilon_{i}^{+},
\]
and
$s^{2}=\operatorname{Var}_{\star}[V_{i}^{+}]$,
equal to that $s^{2}$ in the proof of
Lemma \ref{Lemma_Bootstrapped_Root_Consistency}.

Lemma \ref{Lemma_LOIL_Mean_General} gives $
s^{2}\rightarrow2\sigma^{2},  (P\mbox{-a.s.})
$ and Lemma \ref{Lemma_TCLT_V_Plus_Minus} gives $
\sum_{i=1}^{n_{2}}(V_{i}^{+}-\nu_{i}^{+})/(s\sqrt{n_{2}})\stackrel
{d^{\star}}
{\rightarrow}Z_{1}, (P\mbox{-a.s.})
$. Note that $\sqrt{n}/\sqrt{n_{2}}\rightarrow\sqrt{2/(1-p)}$.
Thus, Slutsky's lemma implies
\[
T_{1}\stackrel{d^{\star}}{\rightarrow}
\frac{\sigma}{f'(d_{0})(1-p)^{1/2}}Z_{1},\qquad
(P\mbox{-a.s.}).
\]
In Lemma \ref{Lemma_Bootstrapped_Root_T2}, following this proof, we
show that for $\gamma\in(0,1/3)$,
\mbox{$
T_{2}\stackrel{P^{\star}}{\rightarrow}0$},  $(P\mbox{-a.s.}).
$
Therefore, another application of Slutsky's lemma completes the proof.
\end{pf*}
\begin{lem}\label{Lemma_Bootstrapped_Root_T2}
For $f\in\mathscr{F}$ and $\gamma\in(0,1/3)$,
$T_{2}\stackrel{P^{\star}}{\rightarrow}0,  (P\mbox{-a.s.})$.
\end{lem}
\begin{pf}
Let
\begin{eqnarray*}
I &=& \hat\beta_{1} - f'(d_{0}),\qquad
\mathit{II} = f(d_{0}) - (2n_{2})^{-1}\sum_{i=1}^{n_{2}}Y_{i}^{+},\\
A &=& \hat\beta_{1}^{\star} -\hat\beta_{1},\qquad
B = (2n_{2})^{-1}\sum_{i=1}^{n_{2}}
(\varepsilon_{i}^{\star+}
- \varepsilon_{i}^{+} ), \\
T_{21} &=& n^{1/2} A\cdot I\cdot \mathit{II},\qquad
T_{22} = n^{1/2} I\cdot B,\\
T_{23} &=& n^{1/2} \mathit{II}\cdot A,\qquad
T_{24} = n^{1/2} A\cdot B.
\end{eqnarray*}
Then
\begin{eqnarray*}
T_{2} &=& n^{1/2}\{-(\hat\beta_{1}^{\star}f'(d_{0}))^{-1}
[ I + A ]\cdot[ \mathit{II} - B ]+ (\hat\beta_{1}f'(d_{0}))^{-1}
I\cdot \mathit{II} \} \\
&=& (\hat\beta_{1}\hat\beta_{1}^{\star}f'(d_{0}))^{-1}n^{1/2}
A\cdot I\cdot \mathit{II}\\
& &{} -
(\hat\beta_{1}^{\star}f'(d_{0}))^{-1}[-n^{1/2} I\cdot B +
n^{1/2} \mathit{II}\cdot A
-n^{1/2} A\cdot B ]\\
&=& (\hat\beta_{1}\hat\beta_{1}^{\star}f'(d_{0}))^{-1}T_{21} -
(\hat\beta_{1}^{\star}f'(d_{0}))^{-1}[-T_{22} + T_{23}
-T_{24} ].
\end{eqnarray*}

It will be shown that
$T_{2i} \stackrel{P^{\star}}{\rightarrow} 0,  (P\mbox{-a.s.})$,
$i=1,2,3,4$ for $\gamma\in(0,1/3)$.
Thus, by Lemma \ref{Lemma_Bootstrapped_Root_Consistency}
and Slutsky' lemma, the conclusion of this lemma holds.

We establish next the convergence of the terms $T_{2i}$.
From
(\ref{TED G beta0 beta1}), (\ref{Bootstrap 2nd stage estimator beta 0 beta1}),
the definitions of $Y'_{i}$, $Y''_{i}$, $Y'^{\star}_{i}$ and
$Y''^{\star}_{i}$, and the Taylor expansions of $f(L)$ and $f(U)$
[(\ref{Two_Ends_Talor_3_prime_1})~and
(\ref{Two_Ends_Talor_3_prime_2})], we have
\begin{eqnarray*}
A
&=&(2Kn_{2})^{-1}n_{1}^{\gamma}\sum_{i=1}^{n_{2}}
( \varepsilon_{i}^{\star-}
-\varepsilon_{i}^{-} )
= (2Kn_{2}^{1/2})^{-1}n_{1}^{\gamma}sn_{2}^{-1/2}\sum_{i=1}^{n_{2}}
( V^{-}_{i} - \nu^{-} )/s,\\
B &=& (2n_{2})^{-1}\sum_{i=1}^{n_{2}}( \varepsilon_{i}^{\star+}
-\varepsilon_{i}^{+} )
=(2n_{2}^{1/2})^{-1}sn_{2}^{-1/2}
\sum_{i=1}^{n_{2}}( V^{+}_{i} - \nu^{+} )/s ,\\
I &=& \hat\beta_{1} - f'(d_{0})
= f''(d_{0})\delta_{d} +
(2K)^{-1}n_{1}^{\gamma}R_{UL}^{-} +
(2Kn_{2})^{-1}n_{1}^{\gamma}\sum_{i=1}^{n_{2}}\varepsilon_{i}^{-},\\
\mathit{II} &=& f(d_{0})-(2n_{2})^{-1}\sum_{i=1}^{n_{2}}Y_{i}^{+} \\
&=& -f'(d_{0})\delta_{d}
-(1/2)f''(d_{0})(\delta_{d}^{2} + K^{2}n_{1}^{-2\gamma} )
-(1/2)R_{UL}^{+}
-(1/n_{2})\sum_{i=1}^{n_{2}}\varepsilon_{i}^{+}.
\end{eqnarray*}
First, consider $T_{21}$. We have
\[
T_{21} = n^{1/2} A\cdot I\cdot \mathit{II} = T'_{21}sn_{2}^{-1/2}\sum_{i=1}^{n_{2}}
( V^{-}_{i} - \nu^{-} )/s,
\]
where $T'_{21} = C_{n}\cdot I\cdot \mathit{II}$ and $C_{n} = n^{1/2}
n_{1}^{\gamma}(2Kn_{2}^{1/2})^{-1}$.
Lemmas \ref{Lemma_LOIL_Mean_General} and
\ref{Lemma_TCLT_V_Plus_Minus} give
\[
s \rightarrow\sqrt{2}\sigma,\qquad (P\mbox{-a.s.}),\qquad
n_{2}^{-1/2}\sum_{i=1}^{n_{2}}
( V^{-}_{i} - \nu^{-} )/s
\stackrel{d^{\star}}{\rightarrow} Z_{2},\qquad (P\mbox{-a.s.}).
\]
Next, it will be shown that $T'_{21}$
converges to 0 $P$-almost surely for $\gamma\in(0,5/12)$.
Then, an application of Slutsky's
lemma gives $T_{21}\stackrel{P^{\star}}{\rightarrow} 0, (P\mbox{-a.s.})$.

With the notation introduced at the beginning of this subsection
and by Lem\-mas~\ref{Lemma LIL one stage estimator} and
\ref{Lemma_LOIL_Mean}, we have, for
$\gamma>0$, $n_{1}^{\gamma}=B_{as}(-\gamma)$,
$(\delta_{d})=B_{as}(1/3)$,
$\sum_{i=1}^{n_{2}}(\varepsilon''_{i}+\varepsilon'_{i})/n_{2}=B_{as}(1/2)$
and
$\sum_{i=1}^{n_{2}}(\varepsilon''_{i}-\varepsilon'_{i})/n_{2}=B_{as}(1/2)$.
Both $R_{U}$ and $R_{L}$ are equal to
$B_{as}(1)+B_{as}(\gamma+ 2/3)+
B_{as}(2\gamma+ 1/3)+B_{as}(3\gamma)$. Thus we have $C_{n} =
B_{as}(-\gamma)$,
$I=B_{as}(1/3)
+ B_{as}(-\gamma)[B_{as}(1)+B_{as}(\gamma+ 2/3)+
B_{as}(2\gamma+ 1/3)+B_{as}(3\gamma)] + B_{as}(1/2-\gamma)=B_{as}(1/3)
+ B_{as}(2\gamma) + B_{as}(1/2-\gamma)$ and
$\mathit{II}= B_{as}(1/3) + [ B_{as}(2/3)+ B_{as}(2\gamma) ]
+ (B_{as}(1)+B_{as}(\gamma+ 2/3)+
B_{as}(2\gamma+ 1/3)+B_{as}(3\gamma)) + B_{as}(1/2)=B_{as}(1/3) +
B_{as}(2\gamma)$.
Thus, for $\gamma\in(0,1/2)$,
\begin{eqnarray*}
T'_{21}
&=& C_{n}\cdot I\cdot \mathit{II}\\
&=& B_{as}(-\gamma)
\times[ B_{as}(1/3)
+ ( B_{as}(2\gamma)) + B_{as}(1/2-\gamma) ]\\
& &{} \times\{ B_{as}(1/3) + B_{as}(2\gamma) \}\\
&=&
B_{as}(2/3-\gamma)+
B_{as}(1/3+\gamma)+
B_{as}(5/6-2\gamma)+
B_{as}(3\gamma).
\end{eqnarray*}
It is easy to see that when $\gamma\in(0,5/12)$, the above upper
bounds $1/2-\gamma$, $1/4+\gamma$, $3/4-2\gamma$, and
$3\gamma$ are all positive. This implies that
$T'_{21}$ converges to 0 $P$-almost surely for $\gamma\in
(0,5/12)$. Therefore, for
$\gamma\in(0,5/12)$, $T_{21}$ converges to 0 in probability
$(P\mbox{-a.s.})$.

Similarly, we can show that $T_{2i}$, $i=2,3$ or $4$, converges to 0
in probability $(P\mbox{-a.s.})$, but with different intervals for
$\gamma$. We next list these results. For $\gamma\in(0,1/2)$,
$T_{22}$ and $T_{24}$ converge to 0 in probability $(P\mbox{-a.s.})$
and for $\gamma\in(0,1/3)$, $T_{23}$ converges to 0 in
probability $(P\mbox{-a.s.})$.
It is worthwhile to note that $\mathscr{F}$ can be considered directly
because the $B_{as}(1/3-\gamma)$ term in $T_{23}$
does not depend on $f''(d_{0})$.
Since $1/3<5/12<1/2$,
$T_{2i}$ converges to 0 in probability $(P\mbox{-a.s.})$ for $i=1,2,3,4$
and $\gamma\in(0,1/3)$.
Thus, for $f\in\mathscr{F}$ and $\gamma\in(0,1/3)$,
$T_{2}$ converges to 0 in probability
$(P\mbox{-a.s.})$.
\end{pf}
\begin{pf*}{Proof of Theorem \ref{Thm Bootstrapped two stage
estimator LD Weak}}
Consider $0 < \gamma< 1/3$. Given an arbitrary subsequence
$\{n_{k}\}_{k=1}^{\infty}$ of $\{n\}_{n=1}^{\infty}$, let $n_{1}=np$ and
$n_{k,1}=n_{k}p$. By Theorem~\ref{Thm_one_stage_IR_estimator_d0}, we
know that $ n_{1}^{\gamma}(\delta_{d})
\equiv(np)^{\gamma}(\hat
d_{np}^{(1)}-d_{0}) \stackrel{P}{\rightarrow} 0$. It follows, by the
relationship between convergence in probability and almost sure
convergence (e.g., see Theorem 20.5 in Billingsley
\cite{Billingsley1995}), that there exists
$\{n_{k(i)}\}_{i=1}^{\infty}$, a further subsequence of $\{n_{k}\}$,
such that $n_{k(i),1}^{\gamma}(\hat d_{n_{k(i),1}}^{(1)}-d_{0})
\rightarrow0, (P\mbox{-a.s.})$.
It now suffices to show that
\[
n_{k(i)}^{1/2}\bigl(\tilde d_{n_{k(i)}}^{(2)\star}-\tilde d_{n_{k(i)}}^{(2)}\bigr)
\stackrel{d^{\star}}{\rightarrow}
C_{2}Z_{1},\qquad
(P\mbox{-a.s.}).
\]

Let $n_{k(i),2}=n_{k(i)}(1-p)/2$. Write $\zeta_{n_{k(i)}}=B_{as}(b)$
if $n_{k(i)}^{\alpha}\zeta_{n_{k(i)}}$ converges to 0 almost surely
for every $\alpha<b$. As in the proof of Theorem \ref{Thm
Bootstrapped two stage estimator LD}, write $n_{k(i)}^{1/2}(\tilde
d_{n_{k(i)}}^{(2)\star}-\tilde d_{n_{k(i)}}^{(2)})$ as $-T_1 + T_2$,
where \textit{both $T_1$ and $T_2$ are now indexed by $n_{k(i)}$}. It
is then not difficult to show that the conditional distribution of
$T_1$ converges to that of $C_{2}Z_{1}$ $P$-almost-surely by
replacing $n$, $n_{1}$ and $n_{2}$ by $n_{k(i)}$, $n_{k(i),1}$ and
$n_{k(i),2}$, respectively, and essentially repeating the steps in Theorem
\ref{Thm Bootstrapped two stage estimator LD}.

It remains to show that $T_{2}\stackrel{P^{\star}}{\rightarrow}0$
$(P\mbox{-a.s.})$. The proof of this follows from that of Lemma
\ref{Lemma_Bootstrapped_Root_T2} by replacing $n$, $n_{1}$ and
$n_{2}$ by $n_{k(i)}$, $n_{k(i),1}$ and $n_{k(i),2}$, respectively,
and noting that $\hat d_{n_{k(i),1}}^{(1)}-d_{0}=B_{as}(1/3)$.
\end{pf*}

%Next, we provide some technical details for Remark \ref{remark for thm
%bootstrpped two stage estimaotor LD}.
%
%4} From the proof of Theorem \ref{Thm Bootstrapped two stage estimator
%LD}
%and Lemma \ref{Lemma_Bootstrapped_Root_T2}, it is seen that
%for $f\in\mathscr{F}$ and $\gamma=1/4$, the slowest term in the
%decomposition of
%$n^{1/2}(\tilde d_{n}^{(2)\star}-\tilde d_{n}^{(2)})$ is
%%
%-\frac{1}{2f'(d_{0})}\frac{\sqrt{n}}{\sqrt{n_{2}}} T^{+}
%+ \frac{1}{2K\hat\beta_{1}^{\star}}\frac{\sqrt{n}}{\sqrt
%{n_{2}}}n_{1}^{1/4}(\hat
%d_{n_{1}}^{(1)}-d_{0}) T^{-},
%%
%where
%%
%T^{+}=\frac{1}{\sqrt{n_{2}}}\sum_{i=1}^{n_{2}}
%[(\varepsilon''^{\star}_{i} + \varepsilon'^{\star}_{i} )-
%(
%%
%and
%%
%T^{-}=\frac{1}{\sqrt{n_{2}}}\sum_{i=1}^{n_{2}} [(\varepsilon
%''^{\star}_{i} -
%)].
%%
%On the other hand, for $\gamma\in(1/4, 1/2)$, the slowest term becomes
%%
%n_{1}^{\gamma}(\hat
%d_{n_{1}}^{(1)}-d_{0}) T^{-}.
%
%By Lemma \ref{Lemma_TCLT_V_Plus_Minus}, we know both $T^{+}$ and
%$T^{-}$ converge weakly to normal distributions $(P-a.s.)$. However,
%it is currently unclear if there exists $\gamma\in[1/4,1/3)$ and a
%real number $c$ (especially 0) such that
%%
%n_{1}^{\gamma}(\hat
%d_{n_{1}}^{(1)}-d_{0})\rightarrow c, (P-a.s.).
%%
%If the answer is positive with $c=0$,
%then Theorem \ref{Thm Bootstrapped two stage estimator LD}
%can be updated by enlarging the upper bound on the range of $\gamma$,
%which will ensure the consistency of the bootstrapped estimator for $f$
%which is not locally linear at $d_{0}$ as well.
%
\vspace{-2pt}
%%%\input{Auxiliary_Lemmas}
%s7.3 ###
\subsection{Some auxiliary lemmas}
First, we state a special almost sure convergence
result on a
triangular array of i.i.d. mean zero random variables.
For the general result,
see Proposition in Hu, M{\'o}ricz and Taylor \cite{Hu1989}.\vspace{-2pt}
\begin{lem}\label{Lemma_LOIL_Mean_General}
If a triangular array of random variables
$\{X_{ni}\}_{i=1}^{m_{n}}$ for $n\in\mathbb{N}$
are i.i.d. copies of a mean $0$ random variable $X$
with $m_{n}$ increases to $\infty$ as $n$ goes to $\infty$
and
$\mathbb{E}|X|^{2p}<\infty$ for some
$p\in[1,2)$,
$
P(\lim_{n\rightarrow\infty}
m_{n}^{-1/p} \sum_{i=1}^{m_{n}}X_{ni} = 0 )= 1$.
\vspace{-2pt}
\end{lem}

Suppose a triangular array of random variables
$\{\varepsilon_{ni}\}_{i=1}^{m_{n}}$ for $n\in\mathbb{N}$
are i.i.d. copies of $\varepsilon$ with mean 0,
where $m_{n}$ increases to $\infty$ as $n$ goes to $\infty$.
Then Lem\-ma~\ref{Lemma_LOIL_Mean_General} tells that
$\bar\varepsilon_{n}=(1/m_{n})\sum_{i=1}^{m_{n}}\varepsilon_{ni}$ and
$(1/m_{n})\sum_{i=1}^{m_{n}}\varepsilon_{ni}^{2}$
converge to $0$ and $\sigma^{2}$ almost surely
given $\mathbb{E}\varepsilon^{2}<\infty$ and $\mathbb{E}\varepsilon
^{4}<\infty$,
respectively. Further,
the following lemma shows that $n^{1/2}$ is an upper boundary rate of the
almost sure convergence of $\bar\varepsilon_{n}$.\vspace{-2pt}
\begin{lem}\label{Lemma_LOIL_Mean}
If $\mathbb{E} \varepsilon^{4}<\infty$,
$
P(\lim_{n\rightarrow\infty}
m_{n}^{\alpha}\bar\varepsilon_{n} = 0 )= 1
$ for each $\alpha<1/2$.\vspace{-2pt}
\end{lem}
\begin{pf}
A direct application of Lemma \ref{Lemma_LOIL_Mean_General}
gives that if $\mathbb{E}|\varepsilon|^{2p}<\infty$ for some
$p\in[1,2)$,
$
P(\lim_{n\rightarrow\infty}
m_{n}^{1-1/p}\bar\varepsilon_{n} = 0 )= 1
$.
On the other hand, $\mathbb{E}\varepsilon^{4}<\infty$ implies that
$\mathbb{E}|\varepsilon|^{2p}<\infty$ for every $p\in[1,2)$. Thus, the
conclusion follows.\vspace{-2pt}
\end{pf}

Suppose $\{\varepsilon'_{i}\}_{i=1}^{n}$,
$\{\varepsilon''_{i}\}_{i=1}^{n}$,
$\{\varepsilon_{i}^{\prime\star}\}_{i=1}^{n}$ and
$\{\varepsilon_{i}^{\prime\prime\star}\}_{i=1}^{n}$
are the second-stage random errors and
the corresponding bootstrapped ones defined in Section \ref
{Bootstrapped Two-Stage Estimator}.
Note that the subscripts of these random variables indicating the
sample size
are suppressed for the simplicity of notation
and that here ``$n$'' is understood as a dummy variable,
not the total sample size.
Recall $V_{i}^{+}=\varepsilon''^{\star}_{i}+\varepsilon'^{\star}_{i}$,
$\nu^{+}=E_{\star}[V_{i}^{+}]$,
$V_{i}^{-}=\varepsilon''^{\star}_{i}-\varepsilon'^{\star}_{i}$ and
$\nu^{-}=E_{\star}[V_{i}^{-}]$, where $E_{\star}$ means the
expectation conditioning on the second-stage data.
Since $\operatorname{Var}_{\star}[V_{i}^{+}]=\operatorname{Var}_{\star}[V_{i}^{-}]$, we
denote both as $s^{2}$.
The following lemma shows that both
$V_{i}^{+}$ and $V_{i}^{-}$ are asymptotically normal $P$-almost
surely.\vspace{-2pt}
\begin{lem}\label{Lemma_TCLT_V_Plus_Minus}
If $\mathbb{E}\varepsilon^{6} <\infty$, we have
\begin{eqnarray*}
\frac{1}{\sqrt{n}}\sum_{i=1}^{n}\frac{V_{i}^{+}-\nu^{+}}{s}
&\stackrel{d^{\star}}{\rightarrow}& Z,\qquad
(P\mbox{-a.s.}),\\[-2pt]
\frac{1}{\sqrt{n}}\sum_{i=1}^{n}\frac{V_{i}^{-}-\nu^{-}}{s}
&\stackrel{d^{\star}}{\rightarrow}& Z,\qquad (P\mbox{-a.s.}),
\end{eqnarray*}
where $Z$ follows a $N(0,1)$ distribution.
\end{lem}
\begin{pf}
We only prove the former
and the latter can be shown similarly.
Let $\xi_{ni}=(V_{i}^{+}-\nu^{+})/(\sqrt{n}s)$, for $i=1,2,\ldots
,n$, and
$S_{n}=\sum_{i=1}^{n}\xi_{ni}$. It is easy to see that $E_{\star
}[\xi_{ni}]=0$
and $\operatorname{Var}_{\star}[S_{n}]=1$. Thus, it suffices to check that the following
Lindeberg condition holds for each $\eta>0$ (see, e.g.,
Theorem 2 on page 334 of Shiryaev \cite{Shiryaev1995}):
$
\sum_{i}^{n}E_{\star}[\xi_{ni}^{2}\{|\xi_{ni}|\geq\eta\}]
\rightarrow0, (P\mbox{-a.s.}).
$
Note that\vspace{-2pt}
\begin{eqnarray*}
\sum_{i}^{n}E_{\star}[\xi_{ni}^{2}\{|\xi_{ni}|\geq\eta\}]
&=& E_{\star}\bigl([(V_{1}^{+}-\nu^{+})/s]^{2}
\bigl\{|(V_{1}^{+}-\nu^{+})/s|\geq\sqrt{n}\eta\bigr\}\bigr)\\[-2pt]
&\leq& \bigl(\sqrt{n}\eta\bigr)^{-1}|s|^{-3}
E_{\star}|V_{1}^{+}-\nu^{+}|^{3},\vspace{-2pt}
\end{eqnarray*}
\[
s^{2}=\frac{1}{n}\sum_{i=1}^{n}(\varepsilon''_{i})^{2}
-\Biggl(\frac{1}{n}\sum_{i=1}^{n}\varepsilon''_{i}\Biggr)^{2}
+\frac{1}{n}\sum_{i=1}^{n}(\varepsilon'_{i})^{2}
-\Biggl(\frac{1}{n}\sum_{i=1}^{n}\varepsilon'_{i}\Biggr)^{2}
\rightarrow
2\sigma^{2}\!,\quad (P\mbox{-a.s.}),
\]
then it is sufficient to show
$\overline{\lim}_{n\rightarrow\infty}
E_{\star}|V_{1}^{+}-\nu^{+}|^{3}<\infty, (P\mbox{-a.s.})$. Since
\begin{eqnarray*}
E_{\star}|V_{1}^{+}-\nu^{+}|^{3}
&\leq& E_{\star} [|V_{1}^{+}|^{3}+|\nu^{+}|^{3} +
3|V_{1}^{+}|^{2}|\nu^{+}| + 3|V_{1}^{+}||\nu^{+}|^{2}]\\
&=& E_{\star} |V_{1}^{+}|^{3}
+3|\nu^{+}|E_{\star} |V_{1}^{+}|^{2}
+3|\nu^{+}|^{2}E_{\star} |V_{1}^{+}|
+|\nu^{+}|^{3},
\end{eqnarray*}
and $\nu^{+}=\frac{1}{n}\sum_{i=1}^{n_{2}}(\varepsilon''_{i} +
\varepsilon'_{i} )\rightarrow0, (P\mbox{-a.s.})$,
it
suffices to show $\overline{\lim}_{n\rightarrow\infty}
E_{\star}|V_{1}^{+}|^{k}<\infty, (P\mbox{-a.s.})$, for $k=1,2,3$.
We only need to show the case where $k=3$.
From $(a+b)^{3}\leq4(a^{3}+b^{3})$ for
nonnegative $a$ and $b$,
\[
E_{\star}|V_{1}^{+}|^{3}
= \frac{1}{n^{2}} \sum_{i=1}^{n} \sum_{j=1}^{n} |\varepsilon
''_{i}+\varepsilon'_{j}|^{3}
\leq4\Biggl(\frac{1}{n}\sum_{i=1}^{n}|\varepsilon''_{i}|^{3}
+ \frac{1}{n}\sum_{i=1}^{n}|\varepsilon'_{i}|^{3} \Biggr).
\]
By Lemma \ref{Lemma_LOIL_Mean_General}, both
$(1/n)\sum_{i=1}^{n}|\varepsilon''_{i}|^{3}$ and
$(1/n)\sum_{i=1}^{n}|\varepsilon'_{i}|^{3}$ converges almost surely
under the assumption $\mathbb{E}\varepsilon^{6}<\infty$.
Therefore, $\overline{\lim}_{n\rightarrow\infty}
E_{\star}|V_{1}^{+}|^{3}<\infty$, $(P\mbox{-a.s.})$,
which completes the proof.
\end{pf}
\end{appendix}

\section*{Acknowledgments}
The authors would like to thank the Editor, the Associate Editor and
two anonymous referees for their constructive comments and suggestions.
In particular, one of the referees pointed out a strategy and a reference
to Durot (2008) that enabled us to strengthen the results on bootstrap
consistency. Further, we would like to acknowledge
Professor Michael Woodroofe's suggestion of
employing a local linear approximation in a shrinking neighborhood
of the target quantity and for many useful discussions on this
topic.

%suskaldyti doi

%
\printaddresses

\end{document}